\documentclass[11pt,leqno]{amsart} 

\usepackage{amscd,fullpage}

\input epic.sty

\newcommand{\Proof}{\begin{proof}}

\newcommand{\Endproof}{\end{proof}}

\renewcommand{\bold}[1]{\textup{\textbf{#1}}}

\newtheorem{theorem}{Theorem}[section]
\newtheorem{proposition}[theorem]{Proposition}
\newtheorem{corollary}[theorem]{Corollary}
\newtheorem{lemma}[theorem]{Lemma}
\newtheorem{definition}[theorem]{Definition}

\newcommand{\R}{\mathbb{R}}
\newcommand{\Q}{\mathbb{Q}}

\newcommand{\p}{\partial}
\newcommand{\half}{\tfrac{1}{2}}

\newcommand{\<}{\langle}
\renewcommand{\>}{\rangle}

\DeclareMathOperator{\ad}{ad}
\DeclareMathOperator{\dga}{\mathsf{dAlg}}
\DeclareMathOperator{\ga}{\mathsf{Alg}}
\let\ss=\undefined
\DeclareMathOperator{\ss}{\mathsf{sSet}}
\DeclareMathOperator{\Chain}{\mathsf{Chain}}
\renewcommand{\o}{\otimes}
\newcommand{\om}{\omega}
\newcommand{\Om}{\Omega}
\newcommand{\bull}{\bullet}

\DeclareMathOperator{\Id}{id}

\renewcommand{\[}{{[\![}}
\DeclareMathOperator{\End}{End}

\newcommand{\eps}{\varepsilon}

\renewcommand{\*}{\cdot}

\newcommand{\Linf}{$L_\infty$}

\newcommand{\g}{\mathfrak{g}}

\newcommand{\h}{\mathfrak{h}}
\newcommand{\m}{\mathbf{m}}
\newcommand{\DELTA}{\mathbf{\Delta}}
\newcommand{\GAMMA}{\gamma}
\newcommand{\K}{\mathsf{K}}
\DeclareMathOperator{\MC}{\mathsf{MC}}
\DeclareMathOperator{\mc}{\mathsf{mc}}

\newcommand{\CC}{\mathcal{C}}

\newcommand{\Wedge}{\Lambda}

\DeclareMathOperator{\e}{\mathbf{e}}
\DeclareMathOperator{\sk}{sk}

\DeclareMathOperator{\cosk}{cosk}
\DeclareMathOperator{\Mor}{Mor}
\DeclareMathOperator{\Map}{Map}
\DeclareMathOperator{\Spec}{Spec}
\newcommand{\op}{\textup{op}}
\renewcommand{\t}{\mathbf{t}}

\newcommand{\CA}{\mathcal{A}}

\newcommand{\ts}{\tilde{s}}
\newcommand{\sbar}{\bar{s}}
\renewcommand{\phi}{\varphi}
\newcommand{\CF}{\mathcal{F}}

\begin{document}

\title{Lie theory for nilpotent \Linf-algebras}

\dedicatory{Dedicated to Ross Street on his sixtieth birthday}

\thanks{This paper was written while I was a guest of R.~Pandharipande
  at Princeton University, with partial support from the Packard
  Foundation. I received further support from the IAS, and from NSF
  Grant DMS-0072508. I am grateful to D. Roytenberg for the suggestion
  to extend Deligne's groupoid to \Linf-algebras.}

\author{Ezra Getzler}

\address{Department of Mathematics, Northwestern University,
  Evanston, Illinois}

\email{getzler@northwestern.edu}

\maketitle

\section{Introduction}

Let $A$ be a differential graded (dg) commutative algebra over a field
$\K$ of characteristic $0$. Let $\Om_\bull$ be the simplicial dg
commutative algebra over $\K$ whose $n$-simplices are the algebraic
differential forms on the $n$-simplex $\Delta^n$. In \cite{Sullivan},
\S\,8, Sullivan introduced a functor
\begin{equation*}
  A \mapsto \Spec_\bull(A) = \dga(A,\Om_\bull)
\end{equation*}
from dg commutative algebras to simplicial sets; here, $\dga(A,B)$ is
the set of morphisms of dg algebras from $A$ to $B$. (Sullivan use the
notation $\<A\>$ for this functor.) This functor generalizes the
spectrum, in the sense that if $A$ is a commutative algebra,
$\Spec_\bull(A)$ is the discrete simplicial set
\begin{equation*}
  \Spec(A)=\ga(A,\K) ,
\end{equation*}
where $\ga(A,B)$ is the set of morphisms of algebras from $A$ to $B$.

If $E$ is a flat vector bundle on a manifold $M$, the complex of
differential forms $\bigl(\Om^*(M,E),d\bigr)$ is a dg module for the
dg Lie algebra $\Om^*(M,\End(E))$; denote the action by $\rho$.  To a
one-form $\alpha\in\Om^1(M,\End(E))$ is associated a covariant
derivative
\begin{equation*}
  \nabla = d + \rho(\alpha) : \Om^*(M,E) \to \Om^{*+1}(M,E) .
\end{equation*}
The equation
\begin{equation*}
  \nabla^2 = \rho\bigl( d \alpha + \half [\alpha,\alpha] \bigr)
\end{equation*}
shows that $\nabla$ is a differential if and only if $\alpha$
satisfies the Maurer-Cartan equation
\begin{equation*}
  d\alpha + \half [\alpha,\alpha] = 0 .
\end{equation*}
This example, and others such as the deformation theory of complex
manifolds of Kodaira and Spencer, motivates the introduction of the
Maurer-Cartan set of a dg Lie algebra $\g$ (Nijenhuis and Richardson
\cite{NR}):
\begin{equation*}
  \MC(\g) = \{ \alpha \in \g^1 \mid \delta \alpha + \half
  [\alpha,\alpha] = 0 \} .
\end{equation*}

There is a close relationship between the Maurer-Cartan set and
Sullivan's functor $\Spec_\bull(A)$, which we now explain. The complex
of Chevalley-Eilenberg cochains $C^*(\g)$ of a dg Lie algebra $\g$ is
a dg commutative algebra whose underlying graded commutative algebra
is the graded symmetric algebra $S(\g[1]^\vee)$; here, $\g[1]$ is the
shifted cochain complex $(\g[1]){}^i=\g^{i+1}$, and $\g[1]^\vee$ is
its dual.

If $\g$ is a dg Lie algebra and $\Om$ is a dg commutative algebra, the
tensor product complex $\g\o\Om$ carries a natural structure of a dg
Lie algebra, with bracket
\begin{equation*}
  [x\o a,y\o b] = (-1)^{|a|\,|y|} [x,y] \, ab .
\end{equation*}
\begin{proposition}
  Let $\g$ be a dg Lie algebra whose underlying cochain complex is
  bounded below and finite-dimensional in each degree. There is a
  natural identification between the $n$-simplices of
  $\Spec_\bull(C^*(\g))$ and the Maurer-Cartan elements of
  $\g\o\Om_n$.
\end{proposition}
\Proof
  Under the stated hypotheses on $\g$, there is a natural
  identification
  \begin{equation*}
    \MC(\g\o\Om) \cong \dga(C^*(\g),\Om)
  \end{equation*}
  for any dg commutative algebra $\Om$. Indeed, there is an inclusion
  \begin{align*}
    \dga(C^*(\g),\Om) &\subset \ga(C^*(\g),\Om) = \ga(S(\g[1]^\vee),\Om) \\
    &\cong ( \g \o \Om )^1 .
  \end{align*}
  It is easily seen that a morphism in $\ga(C^*(\g),\Om)$ is
  compatible with the differentials on $C^*(\g)$ and $\Om$ if and only
  if the corresponding element of $( \g \o \Om )^1$ satisfies the
  Maurer-Cartan equation.
\Endproof

Motivated by this proposition, we introduce for any dg Lie algebra the
simplicial set
\begin{equation*}
  \MC_\bull(\g) = \MC(\g\o\Om_\bull) .
\end{equation*}
According to rational homotopy theory, the functor
$\g\mapsto\MC_\bull(\g)$ induces a correspondence between the homotopy
theories of nilpotent dg Lie algebras over $\Q$ concentrated in
degrees $(-\infty,0]$ and nilpotent rational topological spaces. The
simplicial set $\MC_\bull(\g)$ has been studied in great detail by
Hinich \cite{Hinich}; he calls it the nerve of $\g$ and denotes it by
$\Sigma(\g)$.

However, the simplicial set $\MC_\bull(\g)$ is not the subject of this
paper. Suppose that $\g$ is a nilpotent Lie algebra and let $G$ be the
simply-connected Lie group associated to $\g$. The nerve $N_\bull G$
of $G$ is substantially smaller than $\MC_\bull(\g)$, but they are
homotopy equivalent. In this paper, we construct a natural homotopy
equivalence
\begin{equation}
  \label{embed}
  N_\bull G \hookrightarrow \MC_\bull(\g) ,
\end{equation}
as a special case of a construction applicable to any nilpotent dg Lie
algebra.

To motivate the construction of the embedding \eqref{embed}, we may
start by comparing the sets of $1$-simplices of $N_\bull G$ and of
$\MC_\bull(\g)$. The Maurer-Cartan equation on $\g\o\Om_1$ is
tautologically satisfied, since $\g\o\Om_1$ vanishes in degree $2$;
thus,
\begin{equation*}
  \MC_1(\g) \cong \g[t] dt .  
\end{equation*}
Let $\alpha\in\Om^1(G,\g)$ be the unique left-invariant one-form whose
value $\alpha(e):T_eG\to\g$ at the identity element $e\in G$ is the
natural identification between the tangent space $T_eG$ of $G$ at $e$
and its Lie algebra $\g$. Consider the path space
\begin{equation*}
  P_*G = \{ \tau \in \Mor(\mathbb{A}^1,G) \mid \tau(0)=e \}
\end{equation*}
of algebraic morphisms from the affine line $\mathbb{A}^1$ to $G$.
There is an isomorphism between $P_*G$ and the set $\MC_1(\g)$,
induced by associating to a path $\tau:\mathbb{A}^1\to G$ the one-form
$\tau^*\alpha$.

There is a foliation of $P_*G$, whose leaves are the fibres of the
evaluation map $\tau\mapsto\tau(1)$, and whose leaf-space is $G$.
Under the isomorphism between $P_*G$ and $\MC_1(\g)$, this foliation
is simple to characterize: the tangent space to the leaf containing
$\alpha\in\MC_1(\g)$ is the image under the covariant derivative
\begin{equation*}
  \nabla : \g \o \Om^0_1 \to \g \o \Om^1_1 \cong T_\alpha\MC_1(\g)
\end{equation*}
of the subspace
\begin{equation*}
  \{ x \in \g \o \Om^0_1 \mid x(0) = x(1) = 0 \} .
\end{equation*}
The exponential map
\begin{equation*}
  \exp : \g \to G
\end{equation*}
is a bijection for nilpotent Lie algebras; equivalently, each leaf of
this foliation of $\MC_1(\g)$ contains a unique constant one-form. The
embedding $N_1 G\hookrightarrow \MC_1(\g)$ is the inclusion of the
constant one-forms into $\MC_1(\g)$.

What is a correct analogue in higher dimensions for the condition that
a one-form on $\Delta^1$ is constant? Dupont's explicit proof of the
de~Rham theorem \cite{Dupont}, \cite{DupontBook}, relies on a chain
homotopy $s_\bull:\Om_\bull^*\to\Om_\bull^{*-1}$. This homotopy
induces maps $s_n:\g\o\Om^1_n\to\g\o\Om^0_n$, and we impose the gauge
condition $s_n\alpha=0$, which when $n=1$ is the condition that
$\alpha$ is constant. The main theorem of this paper shows that the
simplicial set
\begin{equation} \label{Gamma}
  \GAMMA_\bull(\g) = \{ \alpha \in \MC_\bull(\g) \mid s_\bull \alpha=0 \}
\end{equation}
is isomorphic to the nerve $N_\bull G$.

The key to the proof of this isomorphism is the verification that
$\GAMMA_\bull(\g)$ is a Kan complex, that is, that it satisfies the
extension condition in all dimensions. In fact, we give explicit
formulas for the required extensions, which yield in particular a new
approach to the Campbell-Hausforff formula.

The definition of $\GAMMA_\bull(\g)$ works \emph{mutatis mutandi} for
nilpotent dg Lie algebras; we argue that $\GAMMA_\bull(\g)$ is a good
generalization to the differential graded setting of the Lie group
associated to a nilpotent Lie algebra. For example, when $\g$ is a
nilpotent dg Lie algebra concentrated in degrees $[0,\infty)$,
$\GAMMA_\bull(\g)$ is isomorphic to the nerve of the Deligne groupoid
$\CC(\g)$.

Recall the definition of this groupoid (cf.\ Goldman and Millson
\cite{GM}). Let $G$ be the nilpotent Lie group associated to the
nilpotent Lie algebra $\g^0\subset\g$.  This Lie group acts on
$\MC(\g)$ by the formula
\begin{equation} \label{Deligne}
  e^X \* \alpha = \alpha - \sum_{n=0}^\infty
  \frac{\ad(X)^n(\delta_\alpha X)}{(n+1)!} .
\end{equation}
The Deligne groupoid $\CC(\g)$ of $\g$ is the groupoid associated to
this group action. There is a natural identification between
$\pi_0(\MC_\bull(\g))$ and $\pi_0(\CC(\g))=\MC(\g)/G$. Following
Kodaira and Spencer, we see that this groupoid may be used to study
the formal deformation theory of such geometric structures as complex
structures on a manifold, holomorphic structures on a complex vector
bundle over a complex manifold, and flat connections on a real vector
bundle.

In all of these cases, the dg Lie algebra $\g$ controlling the
deformation theory is concentrated in degrees $[0,\infty)$, and the
associated formal moduli space is $\pi_0(\MC_\bull(\g))$.  On the
other hand, in the deformation theory of Poisson structures on a
manifold, the associated dg Lie algebra, known as the Schouten Lie
algebra, is concentrated in degrees $[-1,\infty)$.  Thus, the theory
of the Deligne groupoid does not apply, and in fact the formal
deformation theory is modelled by a $2$-groupoid. (This $2$-groupoid
was constructed by Deligne \cite{Deligne}, and, independently, in
Section~2 of Getzler~\cite{darboux}.) The functor $\GAMMA_\bull(\g)$
allows the construction of a candidate Deligne $\ell$-groupoid, if the
nilpotent dg Lie algebra $\g$ is concentrated in degrees
$(-\ell,\infty)$. We present the theory of $\ell$-groupoids in
Section~2, following Duskin \cite{Duskin}, \cite{DuskinNerve} closely.

It seemed most natural in writing this paper to work from the outset
with a generalization of dg Lie algebras called \Linf-algebras. We
recall the definition of \Linf-algebras in Section~4; these are
similar to dg Lie algebras, except that they have a graded
antisymmetric bracket $[x_1,\dotsc,x_k]$, of degree $2-k$, for each
$k$. In the setting of \Linf-algebras, the definition of a
Maurer-Cartan element becomes
\begin{equation*}
  \delta\alpha + \sum_{k=2}^\infty \frac{1}{k!} \,
  [\underbrace{\alpha,\dotsc,\alpha}_{\text{$k$ times}}] = 0 .
\end{equation*}
Given a nilpotent \Linf-algebra $\g$, we define a simplicial set
$\GAMMA_\bull(\g)$, whose $n$-simplices are Maurer-Cartan elements
$\alpha\in\g\o\Om_n$ such that $s_n\alpha=0$. We prove that
$\GAMMA_\bull(\g)$ is a Kan complex, and that the inclusion
$\GAMMA_\bull(\g)\hookrightarrow\MC_\bull(\g)$ is a homotopy
equivalence.

The Dold-Kan functor $K_\bull(V)$ (Dold \cite{Dold}, Kan \cite{Kan})
is a functor from positively graded chain complexes (or equivalently,
negatively graded cochain complexes) to simplicial abelian groups. The
set of $n$-simplices of $K_n(V)$ is the abelian group
\begin{equation}
  \label{kan}
  K_n(V) = \Chain(C_*(\Delta^n),V)
\end{equation}
of morphisms of chain complexes from the complex $C_*(\Delta^n)$ of
normalized simplicial chains on the simplicial set $\Delta^n$ to $V$.
Eilenberg-MacLane spaces are obtained when the chain complex is
concentrated in a single degree (Eilenberg-MacLane \cite{EM}).

The functor $\GAMMA_\bull(\g)$ is a nonabelian analogue of the
Dold-Kan functor $K_\bull(V)$: if $\g$ is an abelian dg Lie algebra
and concentrated in degrees $(-\infty,1]$, there is a natural
isomorphism between $\GAMMA_\bull(\g)$ and $K_\bull(\g[1])$, since
\eqref{kan} has the equivalent form
\begin{equation*}
  K_n(V) = Z^0(C^*(\Delta^n)\o V,d+\delta) ,
\end{equation*}
where $C^*(\Delta^n)$ is the complex of normalized simplicial cochains
on the simplicial set $\Delta^n$.

The functor $\GAMMA_\bull$ has many good features: it carries
surjective morphisms of nilpotent \Linf-algebras to fibrations of
simplicial sets, and carries a large class of weak equivalences of
\Linf-algebras to homotopy equivalences. And of course, it yields
generalizations of the Deligne groupoid, and of the Deligne
$2$-groupoid, for \Linf-algebras. It shares with $\MC_\bull$ an
additional property: there is an action of the symmetric group
$S_{n+1}$ on the set of $n$-simplices $\gamma_n(\g)$ making
$\gamma_\bull$ into a functor from \Linf-algebras to symmetric sets,
in the sense of Fiedorowicz and Loday \cite{FL}. In order to simplify
the discussion, we have not emphasized this point, but this perhaps
indicates that the correct setting for $\ell$-groupoids is the
category of symmetric sets.

\section{Kan complexes and $\ell$-groupoids}

Kan complexes are a natural non-abelian analogue of chain complexes:
just as the homology groups of chain complexes are defined by imposing
an equivalence relation on a subset of the chains, the homotopy groups
of Kan complexes are defined by imposing an equivalence relation on a
subset of the simplices.

Recall the definition of the category of simplicial sets. Let $\Delta$
be the category of finite non-empty totally ordered sets. This
category $\Delta$ has a skeleton whose objects are the ordinals
$[n]=(0<1<\dotsb<n)$; this skeleton is generated by the face maps
$d_k:[n-1]\to[n]$, $0\le k\le n$, which are the injective maps
\begin{equation*}
  d_k(i) = \begin{cases} i , & i<k , \\ i+1 , & i\ge k , \end{cases}
\end{equation*}
and the degeneracy maps $s_k:[n]\to[n-1]$, $0\le k\le n-1$, which are
the surjective maps
\begin{equation*}
  s_k(i) = \begin{cases} i , & i\le k , \\ i-1 , & i>k . \end{cases}
\end{equation*}

A simplicial set $X_\bull$ is a contravariant functor from $\Delta$ to
the category of sets. This amounts to a sequence of sets $X_n=X([n])$
indexed by the natural numbers $n\in\{0,1,2,\dotso\}$, and maps
\begin{align*}
  & \delta_k = X(d_k) : X_n \to X_{n-1} , & 0\le k\le n , \\
  & \sigma_k = X(s_k) : X_{n-1} \to X_n , & 0\le k\le n ,
\end{align*}
satisfying certain relations. (See May \cite{May} for more details.)
A degenerate simplex is one of the form $\sigma_ix$; a nondegenerate
simplex is one which is not degenerate. Simplicial sets form a
category; we denote by $\ss(X_\bull,Y_\bull)$ the set of morphisms
between two simplicial sets $X_\bull$ and $Y_\bull$.

The geometric $n$-simplex $\DELTA^n$ is the convex hull of the unit
vectors $e_k$ in $\R^{n+1}$:
\begin{equation*}
  \DELTA^n = \{ (t_0,\dotsc,t_n) \in [0,1]^{n+1} \mid t_0+\dotsb+t_n=1
  \} .
\end{equation*}
Its $\binom{n+1}{k+1}$ faces of dimension $k$ are the convex hulls of
the nonempty subsets of $\{e_0,\dotsc,e_n\}$ of cardinality $k+1$.

The $n$-simplex $\Delta^n$ is the representable simplicial set
$\Delta^n=\Delta(-,[n])$. Thus, the nondegenerate simplices of
$\Delta^n$ correspond to the faces of the geometric simplex
$\DELTA^n$. By the Yoneda lemma, $\ss(\Delta^n,X_\bull)$ is naturally
isomorphic to $X_n$.

Let $\Delta[k]$ denote the full subcategory of $\Delta$ whose objects
are the simplices $\{[i]\mid i\le k\}$, and let $\sk_k$ be the
restriction of a simplicial set from $\Delta^\op$ to
$\Delta[k]^\op$. The functor $\sk_k$ has a right adjoint $\cosk_k$,
called the $k$-coskeleton, and we have
\begin{equation*}
  \cosk_k(\sk_k(X))_n = \ss(\sk_k(\Delta^n),X_\bull) .
\end{equation*}

For $0\le i\le n$, let $\Lambda_i^n\subset\Delta^n$ be the union of
the faces $d_k[\Delta^{n-1}]\subset\Delta^n$, $k\ne i$. An $n$-horn in
$X_\bull$ is a simplicial map from $\Lambda_i^n$ to $X_\bull$, or
equivalently, a sequence of elements
\begin{equation*}
  (x_0,\dotsc,x_{i-1},-,x_{i+1},\dotsc,x_n)\in(X_{n-1})^n
\end{equation*}
such that $\p_jx_k=\p_{k-1}x_j$ for $0\le j<k\le n$.
\begin{definition}
  A map $f:X_\bull\to Y_\bull$ of simplicial sets is called a
  fibration if the maps
  \begin{equation*}
    \xi^n_i : X_n \to \ss(\Lambda_i^n,X_\bull)
    \times_{\ss(\Lambda^n_i,Y_\bull)} Y_n
  \end{equation*}
  defined by
  \begin{equation*}
    \xi^n_i(x) = (\p_0x,\dotsc,\p_{i-1}x,-,\p_{i+1}x,\dotsc,\p_nx)
    \times f(x)
  \end{equation*}
  are surjective for $n>0$. A simplicial set $X_\bull$ is a Kan
  complex if the map from $X_\bull$ to the terminal object $\Delta^0$
  is a fibration.
\end{definition}

A groupoid is a small category with invertible morphisms. Denote the
sets of objects and morphisms of a groupoid $G$ by $G_0$ and $G_1$,
the source and target maps by $s:G_1\to G_0$ and $t:G_0\to G_1$, and
the identity map by $e:G_0\to G_1$. The nerve $N_\bull G$ of a
groupoid $G$ is the simplicial set whose $0$-simplices are the objects
$G_0$ of $G$, and whose $n$-simplices for $n>0$ are the composable
chains of $n$ morphisms in $G$:
\begin{equation*}
  N_nG=\{[g_1,\dotsc,g_n]\in(G_1)^n\mid sg_i=tg_{i+1}\} .
\end{equation*}
The face and degeneracy maps are defined using the product and the
identity of the groupoid:
\begin{align*}
  \p_k[g_1,\dotsc,g_n] &=
  \begin{cases}
    [g_2,\dotsc,g_n] , & k=0 , \\
    [g_1,\dotsc,g_kg_{k+1},\dotsc,g_n] , & 0<k<n , \\
    [g_1,\dotsc,g_{n-1}] , & k=n ,
  \end{cases} \\
  \sigma_k[g_1,\dotsc,g_{n-1}] &=
  \begin{cases}
    [etg_1,g_1,\dotsc,g_{n-1}] , & k=0 , \\
    [g_1,\dotsc,g_{k-1},etg_k,g_k,\dotsc,g_{n-1}] , & 0<k<n , \\
    [g_1,\dotsc,g_{n-1},esg_{n-1}] , & k=n .
  \end{cases}
\end{align*}

The following characterization of the nerves of groupoids was
discovered by Grothendieck; we sketch the proof.
\begin{proposition} \label{groupoid}
  A simplicial set $X_\bull$ is the nerve of a groupoid if and only if
  the maps
  \begin{equation*}
    \xi^n_i : X_n \to \ss(\Lambda^n_i,X_\bull)
  \end{equation*}
  are bijective for all $n>1$.
\end{proposition}
\Proof
  The nerve of a groupoid is a Kan complex; in fact, it is a very
  special kind of Kan complex, for which the maps $\xi^2_i$ are not
  just surjective, but bijective. The unique filler of the horn
  $(-,g,h)$ is the $2$-simplex $[h,h^{-1}g]$, the unique filler of the
  horn $(g,-,h)$ is the $2$-simplex $[h,g]$, and the unique filler of
  the horn $(g,h,-)$ is the $2$-simplex $[hg^{-1},g]$. Thus, the
  uniqueness of fillers in dimension $2$ exactly captures the
  associativity of the groupoid and the existence of inverses.

  The nerve of a groupoid is determined by its $2$-skeleton, in the
  sense that
  \begin{equation} 
    \label{coskeletal}
    N_\bull G \cong \cosk_2(\sk_2(N_\bull G)) .
  \end{equation}
  It follows from \eqref{coskeletal}, and the bijectivity of the maps
  $\xi^2_i$, that the maps $\xi^n_i$ are bijective for all $n>1$.
  
  Conversely, given a Kan complex $X_\bull$ such that $\xi^n_i$ is
  bijective for $n>1$, we construct a groupoid $G$ such that
  $X_\bull\cong N_\bull G$: $G_i=X_i$ for $i=0,1$, $s=\p_1:G_1\to
  G_0$, $t=\p_0:G_1\to G_0$ and $e=\sigma_0:G_0\to G_1$.
  
  Denote by $\<x_0,\dotsc,x_{k-1},-,x_{k+1},\dotsc,x_n\>$ the unique
  $n$-simplex which fills the horn
  \begin{equation*}
    (x_0,\dotsc,x_{i-1},-,x_{i+1},\dotsc,x_n)\in\ss(\Lambda^n_i,X) .
  \end{equation*}
  Given a pair of morphisms $g_1,g_2\in G_1$ such that $sg_1=tg_2$,
  define their composition by the formula
  \begin{equation*}
    g_1g_2 = \p_1\<g_2,-,g_1\> .
  \end{equation*}
  Given three morphisms $g_1,g_2,g_3\in G_1$ such that $sg_1=tg_2$ and
  $sg_2=tg_3$, the $3$-simplex $x = [g_1,g_2,g_3] \in X_3$ satisfies
  \begin{equation*}
    g_1(g_2g_3) = \p_1\p_2x = \p_1\p_1x = (g_1g_2)g_3 ,
  \end{equation*}
  hence composition in $G_1$ is associative. For a picture of the
  $3$-simplex $x$, see Figure~1.

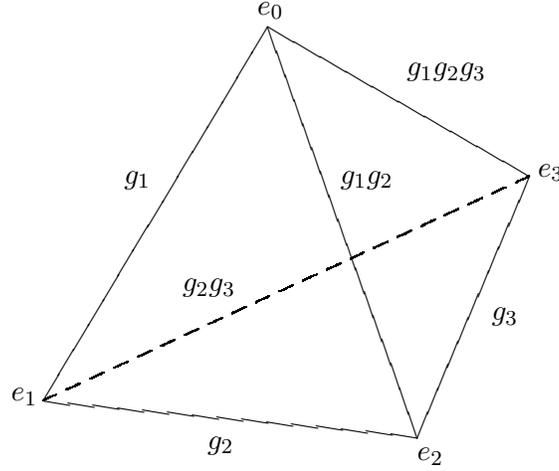
\begin{figure}[h]
  \centering
  \setlength{\unitlength}{0.00087489in}
{\renewcommand{\dashlinestretch}{30}
\begin{picture}(2949,2803)(0,-10)
\put(500,1800){\makebox(0,0)[lb]{$g_1$}}
\put(1000,200){\makebox(0,0)[lb]{$g_2$}}
\put(2712,976){\makebox(0,0)[lb]{$g_3$}}
\put(1800,1800){\makebox(0,0)[lb]{$g_1g_2$}}
\put(850,1150){\makebox(0,0)[lb]{$g_2g_3$}}
\put(2200,2450){\makebox(0,0)[lb]{$g_1g_2g_3$}}
\put(1300,2825){\makebox(0,0)[lb]{$e_0$}}
\put(-180,500){\makebox(0,0)[lb]{$e_1$}}
\put(2262,140){\makebox(0,0)[lb]{$e_2$}}
\put(2990,1850){\makebox(0,0)[lb]{$e_3$}}
\drawline(2262,301)(2937,1876)
\drawline(1362,2776)(2937,1876)
\drawline(1362,2776)(2262,301)
\drawline(12,526)(2262,301)
\drawline(1362,2776)(12,526)
\dashline{96.000}(12,526)(2937,1876)
\end{picture}
}
  \caption{The $3$-simplex $[g_1,g_2,g_3]$}
\end{figure}

  The inverse of a morphism $g\in G_1$ is defined by the formulas
  \begin{equation*}
    g^{-1} = \p_0\<-,etg,g\> = \p_2\<g,esg,-\> .
  \end{equation*}
  To see that these two expressions are equal, call them respectively
  $g^{-\ell}$ and $g^{-\rho}$, and use associativity:
  \begin{equation*}
    g^{-\ell} = g^{-\ell}(gg^{-\rho}) = (g^{-\ell}g)g^{-\rho} =
    g^{-\rho} .
  \end{equation*}
  It follows easily that $(g^{-1})^{-1}=g$, that $g^{-1}g=esg$ and
  $gg^{-1}=etg$, and that $sg^{-1}=tg$ and $tg^{-1}=sg$.
  
  It is clear that
  \begin{equation*}
    sh=\p_1\p_2\<g,-,h\> = \p_1\p_1\<g,-,h\>=s(gh)
  \end{equation*}
  and that
  \begin{equation*}
    tg = \p_0\p_0\<g,-,h\> = \p_0\p_1\<g,-,h\> = t(gh) .
  \end{equation*}
  We also see that
  \begin{align*}
    g &= \p_1\sigma_1[g] = \p_1[g,esg] = g(esg) \\ 
    &= \p_1\sigma_0[g] = \p_1[etg,g] = (etg)g .
  \end{align*}
  Thus, $G$ is a groupoid. Since $\sk_2(X_\bull)\cong\sk_2(N_\bull
  G)$, we conclude by \eqref{coskeletal} that $X_\bull\cong N_\bull
  G$.
\Endproof

Duskin has defined a sequence of functors $\Pi_\ell$ from the category
of Kan complexes to itself, which give a functorial realization of the
Postnikov tower. (See Duskin \cite{Duskin} and Glenn \cite{Glenn}, and
for a more extended discussion, Beke \cite{Beke}.)  Let $\sim_\ell$ be
the equivalence relation of homotopy relative to the boundary on the
set $X_\ell$ of $\ell$-simplices. Then $\sk_\ell(X_\bull)/\sim_\ell$
is a well-defined $\ell$-truncated simplicial set, and there is a map
of truncated simplicial sets
\begin{equation*}
  \sk_\ell(X_\bull) \to \sk_\ell(X_\bull)/\sim_\ell ,
\end{equation*}
and by adjunction, a map of simplicial sets
\begin{equation*}
  X_\bull \to \cosk_\ell(\sk_\ell(X_\bull)/\sim_\ell) .  
\end{equation*}
Define $\Pi_\ell(X_\bull)$ to be the image of this map. Then the
functor $\Pi_\ell$ is an idempotent monad on the category of Kan
complexes. If $x_0\in X_0$, we have
\begin{equation*}
  \pi_i(X_\bull,x_0) = \begin{cases}
    \pi_i(\Pi_\ell(X_\bull),x_0) , & i\le\ell , \\
    0 , & i>\ell .
  \end{cases}
\end{equation*}
Thus $\Pi_\ell(X_\bull)$ is a realization of the Postnikov
$\ell$-section of the simplicial set $X_\bull$. For example,
$\Pi_0(X_\bull)$ is the discrete simplicial set $\pi_0(X_\bull)$ and
$\Pi_1(X_\bull)$ is the nerve of the fundamental groupoid of
$X_\bull$. It is interesting to compare $\Pi_\ell(X_\bull)$ to other
realizations of the Postnikov tower, such as
$\cosk_{\ell+1}\bigl(\sk_{\ell+1}\bigl(X_\bull\bigr)\bigr)$: it is a
more economic realization of this homotopy type, and has a more
geometric character.

We now recall Duskin's notion of higher groupoid: he calls these
$\ell$-di\-men\-si\-on\-al hypergroupoids, but we call them simply
weak $\ell$-groupoids.
\begin{definition}
  A Kan complex $X_\bull$ is a weak $\ell$-groupoid if
  $\Pi_\ell(X_\bull)=X_\bull$, or equivalently, if the maps $\xi^n_i$
  are bijective for $n>\ell$; it is a weak $\ell$-group if in addition
  it is reduced (has a single $0$-simplex).
\end{definition}

The $0$-simplices of an $\ell$-groupoid are interpreted as its objects
and the $1$-simplices as its morphisms. The composition $gh$ of a pair
of $1$-morphisms with $\p_1g=\p_0h$ equals $\p_1z$, where $z\in X_2$
is a filler of the horn
\begin{equation*}
  (g,-,h) \in \ss(\Lambda^2_1,X_\bull) .
\end{equation*}
If $\ell>1$, this composition is not canonical --- it depends on the
choice of the filler $z\in X_2$ --- but it is associative up to a
homotopy, by the existence of fillers in dimension $3$.

A weak $0$-groupoid is a discrete set, while a weak $1$-groupoid is
the nerve of a groupoid, by Proposition \ref{groupoid}. In
\cite{DuskinNerve}, Duskin identifies weak $2$-groupoids with the
nerves of bigroupoids. A bigroupoid $G$ is a bicategory whose
$2$-morphisms are invertible, and whose $1$-morphisms are
equivalences; the nerve $N_\bull G$ of $G$ is a simplicial set whose
$0$-simplices are the objects of $G$, whose $1$-simplices are the
morphism of $G$, and whose $2$-simplices $x$ are the $2$-morphisms
with source $\p_2x\circ\p_0x$ and target $\p_1x$.

The singular complex of a topological space is the simplicial set
\begin{equation*}
  S_n(X) = \Map(\DELTA^n,X) .
\end{equation*}
To see that this is a Kan complex, we observe that there is a
continuous retraction from $\DELTA^n=|\Delta^n|$ to $|\Lambda^n_i|$.
The fundamental $\ell$-groupoid of a topological space $X$ is the weak
$\ell$-groupoid $\Pi_\ell(S_\bull(X))$. For $\ell=0$, this equals
$\pi_0(X)$, while for $\ell=1$, it is the nerve of the fundamental
groupoid of $X$.

Often, weak $\ell$-groupoids come with explicit choices for fillers of
horns: tentatively, we refer to such weak $\ell$-groupoids as
$\ell$-groupoids. (Often, this term is used for what we call strict
$\ell$-groupoids, but the latter are of little interest for $\ell>2$.)
We may axiomatize $\ell$-groupoids by a weakened form of the axioms
for simplicial $T$-complexes, studied by Dakin \cite{Dakin} and Ashley
\cite{Ashley}.
\begin{definition}
  An \bold{$\ell$-groupoid} is a simplicial set $X_\bull$ together
  with a set of thin elements $T_n\subset X_n$ for each $n>0$,
  satisfying the following conditions:
  \begin{enumerate}
  \item every degenerate simplex is thin;
  \item every horn has a unique thin filler;
  \item every $n$-simplex is thin if $n>\ell$.
  \end{enumerate}
\end{definition}

If $\g$ is an $\ell$-groupoid and $n>\ell$, we denote by
$\<x_0,\dotsc,x_{i-1},-,x_{i+1},\dotsc,x_n\>$ the unique thin filler of
the horn
\begin{equation*}
  (x_0,\dotsc,x_{i-1},-,x_{i+1},\dotsc,x_n) \in \ss(\Lambda^n_i,X_\bull) .
\end{equation*}

\begin{definition}
  An \bold{$\infty$-groupoid} is a simplicial set $X_\bull$ together
  with a set of thin elements $T_n\subset X_n$ for each $n>0$,
  satisfying the following conditions:
  \begin{enumerate}
  \item every degenerate simplex is thin;
  \item every horn has a unique thin filler.
  \end{enumerate}
\end{definition}

It is clear that every $\ell$-groupoid is a weak $\ell$-groupoid, and
every $\infty$-groupoid is a Kan complex.

\section{The simplicial de~Rham theorem}

Let $\Om_n$ be the free graded commutative algebra over $\K$ with
generators $t_i$ of degree $0$ and $dt_i$ of degree $1$, and relations
$T_n=0$ and $dT_n=0$, where $T_n=t_0+\dotsb+t_n-1$:
\begin{equation*}
  \Om_n = \K[t_0,\dotsc,t_n,dt_0,\dotsc,dt_n] / (T_n,dT_n) .
\end{equation*}
There is a unique differential on $\Om_n$ such that $d(t_i)=dt_i$, and
$d(dt_i)=0$.

The dg commutative algebras $\Om_n$ are the components of a simplicial
dg commutative algebra $\Om_\bull$: the simplicial map $f:[k]\to[n]$
acts by the formula
\begin{equation*}
  f^*t_i = \sum_{f(j)=i} t_j , \quad 0\le i\le n .
\end{equation*}
Using the simplicial dg commutative algebra $\Om_\bull$, we can define
the dg commutative algebra of piecewise polynomial differential forms
$\Om(X_\bull)$ on a simplicial set $X_\bull$. (See Sullivan
\cite{Sullivan}, Bousfield and Guggenheim \cite{BG}, and Dupont
\cite{Dupont}, \cite{DupontBook}.)
\begin{definition}
  The complex of differential forms $\Om(X_\bull)$ on a simplicial set
  $X_\bull$ is the space
  \begin{equation*}
    \Om(X_\bull) = \ss(X_\bull,\Om_\bull)
  \end{equation*}
  of simplicial maps from $X_\bull$ to $\Om_\bull$.
\end{definition}

When $\K=\R$ is the field of real numbers, $\Om(X_\bull)$ may be
identified with the complex of differential forms on the realization
$|X_\bull|$ which are polynomial on each geometric simplex of
$|X_\bull|$.

The following lemma may be found in Bousfield and Guggenheim
\cite{BG}; we learned this short proof from a referee.
\begin{lemma} \label{contractible}
  For each $k\ge0$, the simplicial abelian group $\Om^k_\bull$ is
  contractible.
\end{lemma}
\Proof
  The homotopy groups of the simplicial set $\Om^k_\bull$ equal the
  homology groups of the complex $C_\bull=\Om^k_\bull$ with
  differential
  \begin{equation*}
    \p = \sum_{i=0}^n (-1)^i \p_i : C_n \to C_{n-1} .
  \end{equation*}
  Thus, to prove the lemma, it suffices to construct a contracting
  chain homotopy for the complex $C_\bull$.

  For $0\le i\le n$, let $\pi_i:\DELTA^{n+1}\to\DELTA^n$ be the affine
  map
  \begin{equation*}
    \pi_i(t_0,\dotsc,t_{n+1}) =
    (t_0,\dotsc,t_{i-1},t_i+t_{n+1},t_{i+1},\dotsc,t_n) .
  \end{equation*}
  Define a chain homotopy $\eta:C_n\to C_{n+1}$ by
  \begin{equation*}
    \eta\om = (-1)^{n+1} \sum_{i=0}^n t_i \, \pi_i^*\om .
  \end{equation*}
  For $\om\in\Om^k_n$, we see that
  \begin{equation*}
    \p_i \eta \om =
    \begin{cases}
      - \eta \p_i\om , & 0\le i\le n , \\
      (-1)^{n+1} \, \om , & i=n+1 .
    \end{cases}
  \end{equation*}
  It follows that $(\p \eta + \eta \p)\om=\om$.
\Endproof

Given a sequence $(i_0,\dotsc,i_k)$ of elements of the set
$\{0,\dotsc,n\}$, let
\begin{equation*}
  I_{i_0\dotso i_k} : \Om_n \to \K
\end{equation*}
be the integral over the $k$-chain on the $n$-simplex spanned by the
sequence of vertices $(e_{i_0},\dotsc,e_{i_k})$; this is defined by the
explicit formula
\begin{equation*}
  I_{i_0\dotso i_k}\bigl( t_{i_1}^{a_1} \dotso t_{i_k}^{a_k} \, dt_{i_1}
  \dotsm dt_{i_k} \bigr) = \frac{a_1!\dotsm a_k!}{(a_1+\dotsb+a_k+k)!} .
\end{equation*}
Specializing $\K$ to the field of real numbers, this becomes the usual
Riemann integral.

The space $C_n$ of elementary forms is spanned by the differential
forms
\begin{align*} \label{elementary}
  \om_{i_0\dotso i_k} = k! \, \sum_{j=0}^k (-1)^j t_{i_j} dt_{i_0}
  \dotsm \widehat{dt}_{i_j} \dotsm dt_{i_k} .
\end{align*}
(The coefficient $k!$ normalizes the form in such a way that
$I_{i_0\dotso i_k}(\om_{i_0\dotso i_k})=1$.) The spaces $C_n$ are closed
under the action of the exterior differential,
\begin{equation*}
  d\om_{i_0\dotso i_k} = \sum_{i=0}^n \om_{ii_0\dotso i_k} ,
\end{equation*}
and assemble to a simplicial subcomplex of $\Om_\bull$. The complex
$C_n$ is isomorphic to the complex of simplicial chains on $\Delta^n$,
and this isomorphism is compatible with the simplicial structure. In
\cite{Whitney}, Whitney constructs an explicit projection $P_n$ from
$\Om_n$ to $C_n$:
\begin{equation}
  \label{Pn}
  P_n\om = \sum_{k=0}^n \sum_{i_0<\dotsb<i_k} \om_{i_0\dotso i_k} \,
  I_{i_0\dotso i_k}(\om) .
\end{equation}
The projections $P_n$ assemble to form a morphism of simplicial
cochain complexes $P_\bull:\Om_\bull\to C_\bull$. If $X_\bull$ is a
simplicial set, the complex of elementary forms
\begin{equation*}
  C(X_\bull) = \ss(X_\bull,C_\bull) \subset \Om(X_\bull)
\end{equation*}
on $X_\bull$ is naturally isomorphic to the complex of normalized
simplicial cochains.

\begin{definition}
  A \bold{contraction} is a simplicial endomorphism
  $s_\bull:\Om_\bull^*\to\Om_\bull^{*-1}$ such that
  \begin{equation}
    \label{Dupont}
    \Id - P_\bull = ds_\bull + s_\bull d .
  \end{equation}
\end{definition}

If $X_\bull$ is a simplicial complex, a contraction $s_\bull$ induces
a chain homotopy $s:\Om^*(X_\bull)\to\Om^{*-1}(X_\bull)$ between the
complex of differential forms on $X_\bull$ and the complex
$C(X_\bull)$ of simplicial cochains. In other words, a contraction is
an explicit form of the de~Rham theorem.

Next, we derive some simple properties of a contraction which we will
need later. If $a$ and $b$ are operators on a chain complex
homogeneous of degree $k$ and $\ell$ respectively, we denote by
$[a,b]$ the graded commutator
\begin{equation*}
  [a,b] = ab - (-1)^{k\ell} ba .
\end{equation*}
In particular, of $a$ is homogeneous of odd degree, then
$\half[a,a]=a^2$.
\begin{lemma}
  \label{contraction}
  Let $s_\bull$ be a contraction. Then
  \begin{enumerate}
  \item $P_\bull\,s_\bull=0$;
  \item $s_\bull P_\bull=[d,(s_\bull){}^2]$.
  \end{enumerate}
\end{lemma}
\Proof
  To show that $P_\bull\,s_\bull=0$, we must check that $I_{i_0\dots
    i_k}\circ s_n=0$ for each sequence $(i_0\dotso i_k)$.  By the
  compatibility of $s_\bull$ with simplicial maps, this follows from
  the formula
  \begin{equation*}
    I_{0\dotso k}\circ s_k=0 ,
  \end{equation*}
  which is clear, since $s_k\om$ is a differential form on $\DELTA^k$ of
  degree less than $k$.

  The second part of the lemma is a simple calculation.
\Endproof

Dupont \cite{Dupont}, \cite{DupontBook}, found an explicit
contraction: we now recall his formula. Given $0\le i\le n$, define
the dilation map
\begin{equation*}
  \phi_i : [0,1]\times \DELTA^n \to \DELTA^n
\end{equation*}
by the formula
\begin{equation*}
  \phi_i(u,\t) = u\t + (1-u) e_i .
\end{equation*}
Let $\pi_*:\Om^*([0,1]\times\DELTA^n)\to\Om^{*-1}(\DELTA^n)$ be
integration along the fibers of the projection
$\pi:[0,1]\times\DELTA^n\to\DELTA^n$. Define the operator
\begin{equation*}
  h^i_n : \Om^*_n \to \Om^{*-1}_n
\end{equation*}
by the formula
\begin{equation}
  \label{h}
  h^i_n\om = \pi_*\phi_i^*\om ,
\end{equation}
Let $\eps^i_n:\Om_n\to\K$ be evaluation at the vertex $e_i$. Stokes's
theorem implies the Poincar\'e lemma, that $h^i_n$ is a chain homotopy
between the identity and $\eps^i_n$:
\begin{equation} \label{Poincare}
  dh^i_n + h^i_nd = \Id_n - \eps^i_n .
\end{equation}

The flow $\phi_i(u)$ is generated by the vector field
\begin{equation*}
  E_i = \sum_{j=0}^n (t_j-\delta_{ij}) \p_j .
\end{equation*}
Let $\iota_i$ be the contraction $\iota(E_i)$: we have
\begin{equation}
  \label{iota}
  \iota_j \phi_i(u) = \phi_i(u) \bigl( u \iota_j + (1-u) \iota_i ,
\end{equation}
and also
\begin{equation}
  \label{iom}
  \iota_i \om_{i_0\dotso i_k} = k \, \sum_{p=0}^k (-1)^{p-1} \,
  \delta_{ii_p} \, \om_{i_0\dotso\widehat{\imath}{}_p\dotso i_k} .
\end{equation}
The formula \eqref{h} for $h^i_n$ may be written more explicitly as
\begin{equation*}
  h^i_n = \int_0^1 u^{-1} \, \phi_i(u) \, \iota_i \, du .
\end{equation*}

\begin{lemma}
  \label{h2}
  $h^ih^j+h^jh^i=0$
\end{lemma}
\Proof
  Let $\phi_{ij}:[0,1]\times[0,1]\times\DELTA^n\to\DELTA^n$ be the map
  \begin{equation*}
    \phi_{ij}(u,v,\t) = uvt_k + (1-u)e_i + u(1-v)e_j .
  \end{equation*}
  Then we have
  \begin{equation*}
    h^ih^j\om = \pi_*\phi_{ij}^*\om .
  \end{equation*}
  We have
  \begin{equation*}
    \phi_{ji}(u,v) = \phi_{ij}(\tilde{v},\tilde{u}) ,
  \end{equation*}
  where $\tilde{u}$ and $\tilde{v}$ are determined implicitly by the
  equations
  \begin{align*}
    (1-u)v &= 1 - \tilde{u} , & 1 - v &= (1-\tilde{v})\tilde{u} .
  \end{align*}
  Since this change of variables is a diffeomorphism of the interior
  of the square $[0,1]\times[0,1]$, the lemma follows.
\Endproof

\begin{lemma}
  \label{I}
  $I_{i_0\dotso i_k}(\om) = (-1)^k \, \eps_n^{i_k}h_n^{i_{k-1}} \dots
  h_n^{i_0}\om$
\end{lemma}
\Proof
  For $k=0$, this holds by definition. We argue by induction on $k$.
  We may assume that $\om$ has positive degree, and hence that
  $\om=d\nu$ is exact. By Stokes's theorem,
  \begin{equation*}
    I_{i_0\dotso i_k}(d\nu) = \sum_{j=0}^k (-1)^{j-1}
    I_{i_0\dotso\widehat{\imath}_j\dotso i_k}(\nu) .
  \end{equation*}
  On the other hand, by \eqref{Poincare}, we have
  \begin{align*}
    \eps_n^{i_k} h_n^{i_{k-1}}\dotsm h_n^{i_0} d\nu &=
    \sum_{j=0}^{k-1} (-1)^j \eps_n^{i_k} h_n^{i_{k-1}}\dotsm
    [d,h_n^{i_j}] \dotsm h_n^{i_0} \nu \\
    &= \sum_{j=0}^{k-1} (-1)^j \eps_n^{i_k} h_n^{i_{k-1}}\dotsm
    \widehat{h}{}_n^{i_j} \dotsm h_n^{i_0} \nu + (-1)^k \eps_n^{i_k}
    \eps_n^{i_{k-1}} h_n^{i_{k-2}} \dotsm h_n^{i_0} \nu .
  \end{align*}
  But $\eps_n^{i_k}\eps_n^{i_{k-1}}=\eps_n^{i_{k-1}}$.
\Endproof

The following theorem is due to Dupont.
\begin{theorem}
  The operators
  \begin{equation}
    \label{s}
    s_n = \sum_{k=0}^{n-1} \sum_{i_0<\dots<i_k} \om_{i_0\dotso i_k}
    h_n^{i_k}\dotso h_n^{i_0} , \quad n\ge0 ,
  \end{equation}
  form a contraction.
\end{theorem}
\Proof
  It is straightforward to check that $s_\bull$ is simplicial. In the
  proof of \eqref{Dupont}, we abbreviate $h^i_n$ to $h^i$. In the
  definition of $s_n$, we may take the upper limit of the sum over $k$
  to be $n$. We now have
  \begin{align}
    \label{ds}
    [d,s_n] &= \sum_{k=0}^{n-1} \sum_{i_0<\dots<i_k}
    \sum_{i\notin\{i_0,\dotsc,i_k\}}
    \om_{ii_0\dotso i_k} h^{i_k}\dotsm h^{i_0} \\
    &+ \sum_{k=0}^n \sum_{j=0}^k (-1)^j \sum_{i_0<\dots<i_k}
    \om_{i_0\dotso i_k} h^{i_k}\dotsm [d,h^{i_j}] \dotsm h^{i_0} . \notag
  \end{align}
  By \eqref{Poincare}, we have
  \begin{align*}
    \sum_{k=0}^n \sum_{j=0}^k (-1)^j \sum_{i_0<\dots<i_k}
    \om_{i_0\dotso i_k} h^{i_k}\dotsm [d,h^{i_j}] \dotsm h^{i_0} = \Id
    &+ \sum_{k=1}^n \sum_{j=0}^k (-1)^j \sum_{i_0<\dots<i_k}
    \om_{i_0\dotso i_k} h^{i_k}\dotsm \widehat{h}^{i_j} \dotsm h^{i_0} \\
    &- \sum_{k=0}^n (-1)^k \sum_{i_0<\dots<i_k} \om_{i_0\dotso i_k}
    \eps^{i_k} h^{i_{k-1}} \dotsm h^{i_0} .
  \end{align*}
  The first term on the right-hand side equals the identity operator,
  the second cancels the first sum of \eqref{ds}, while by
  Lemma~\ref{I}, the third sum equals $P_n$.  \Endproof

We will need special class of contractions, which we call gauges.
\begin{definition}
  A \bold{gauge} is a contraction such that $(s_\bull)^2=0$.  
\end{definition}

In fact, Dupont's operator $s_\bull$ is a gauge. But by a trick of
Lambe and Stasheff \cite{LS}, any contraction gives rise to a gauge.
\begin{proposition}
  If $s_\bull$ is a contraction, then the operator
  \begin{equation*}
    \ts_\bull = s_\bull d s_\bull \, (\Id - P_\bull)
  \end{equation*}
  is a gauge. If $s_\bull$ is a gauge, then $\ts_\bull=s_\bull$.
\end{proposition}
\Proof
  Let $\sbar_\bull$ be the contraction
  \begin{equation*}
    \sbar_\bull = s_\bull(\Id-P_\bull) .
  \end{equation*}
  By construction, $\sbar_\bull P_\bull=0$, hence by Lemma
  \ref{contraction}, $[d,(\sbar_\bull){}^2]=0$. Then
  $\ts_\bull=\sbar_\bull d \sbar_\bull$ is a contraction:
  \begin{align*}
    [d,\ts_\bull] &= [d,\sbar_\bull d\sbar_\bull] \\
    &= [d,\sbar_\bull]d\sbar_\bull + \sbar_\bull d [d,\sbar_\bull] \\
    &= (\Id-P_\bull)d\sbar_\bull + \sbar_\bull d (\Id-P_\bull) \\
    &= d(\Id-P_\bull)\sbar_\bull + \sbar_\bull(\Id-P_\bull)d \\
    &= [d,\sbar_\bull] = \Id-P_\bull .
  \end{align*}
  Since $d(\sbar_\bull){}^2d=(\sbar_\bull){}^2d^2=0$, $\ts_\bull$ is a
  gauge:
  \begin{equation*}
    (\ts_\bull){}^2 = (\sbar_\bull d \sbar_\bull)(\sbar_\bull d
    \sbar_\bull) = \sbar_\bull d (\sbar_\bull){}^2 d \sbar_\bull = 0 .
  \end{equation*}

  If $s_\bull$ happens to be a gauge, then by Lemma~\ref{contraction},
  $s_\bull P_\bull=0$. It follows that
  \begin{align*}
    \ts_\bull - s_\bull &= s_\bull ( d s_\bull ( \Id - P_\bull ) - \Id
    ) \\
    &= s_\bull ( d s_\bull - \Id ) \\
    &= - s_\bull ( s_\bull d + P_\bull ) = - (s_\bull){}^2 d + s_\bull
    P_\bull = 0 ,
  \end{align*}
  showing that $\ts_\bull=s_\bull$.
\Endproof

We now turn to the proof that Dupont's operator $s_\bull$ is a gauge.
Denote by $\eps(\alpha)$ the operation of multiplication by a
differential form $\alpha$ on $\Om_n$.
\begin{lemma}
  If $i\notin\{i_0,\dotsc,i_k\}$, then
  \begin{equation*}
    \eps(\om_{i_0\dotso i_k}) h^i = (-1)^k h^i \bigl(
    \eps(\om_{i_0\dotso i_k}) + \eps(\om_{i_0\dotso i_k i}) h^i \bigr) .
  \end{equation*}
\end{lemma}
\Proof
  We have
  \begin{align*}
    (-1)^k h^i \eps(\om_{i_0\dotso i_k}) &= (-1)^k \, \int_0^1 w^{-1}
    \phi_i(w) \iota_i \, \eps(\om_{i_0\dotso i_k}) \, dw \\
    &= \eps(\om_{i_0\dotso i_k}) \int_0^1 w^k \, \phi_i(w) \, \iota_i
    \, dw .
  \end{align*}
  On the other hand, by \eqref{iom},
  \begin{align*}
    (-1)^k h^i \eps(\om_{i_0\dotso i_ki}) h^i &= (-1)^k \int_0^1
    \int_0^1 (uv)^{-1} \phi_i(u) \, \iota_i \, \eps(\om_{i_0\dots
      i_ki}) \, \phi_i(v) \, \iota_i \, dv \, du \\
    &= (k+1) \int_0^1 \int_0^1 (uv)^{-1} \phi_i(u) \,
    \eps(\om_{i_0\dotso i_k}) \, \phi_i(v) \, \iota_i \, dv \,du \\
    &= (k+1) \, \eps(\om_{i_0\dotso i_k}) \int_0^1 \int_0^1 u^k v^{-1}
    \phi_i(uv) \, \iota_i \, dv \,du .
  \end{align*}
  Changing variables from $u$ to $w=uv$, we see that
  \begin{align*}
    \int_0^1 \int_0^1 u^k v^{-1} \phi_i(uv) \, dv \,du &= \int_0^1
    \biggl( \int_w^1 v^{-k-2} \, dv \biggr) w^k \, \phi_i(w) \,dw \\
    &= (k+1)^{-1} \int_0^1 ( w^{-1} - w^k ) \phi_i(w) \,dw ,
  \end{align*}
  establishing the lemma.
\Endproof

\begin{theorem}
  \label{gauge}
  The operator $s_\bull$ is a gauge.
\end{theorem}
\Proof
  By induction on $k$, the above lemma shows that
  \begin{equation*}
    h^{i_k} \dotsm h^{i_0} s = \sum_{\ell=0}^{n-1} (-1)^{k\ell+\ell}
    \sum_{\substack{j_0<\dots<j_\ell \\
        \{i_0,\dotsc,i_k\}\cap\{j_0,\dotsc,j_\ell\}=\emptyset}}
    \om_{j_0\dotso j_\ell} h^{i_k} \dotsm h^{i_0} h^{j_\ell} \dotsm h^{j_0} .
  \end{equation*}
  It follows that $s^2$ is given by the formula
  \begin{equation}
    \label{s2}
    s^2 = \sum_{k,\ell=0}^\infty (-1)^{k\ell+\ell}
    \sum_{\substack{i_0<\dots<i_k ; j_0<\dots<j_\ell \\
        \{i_0,\dotsc,i_k\}\cap\{j_0,\dotsc,j_\ell\}=\emptyset}}
    \om_{i_0\dotso i_k} \om_{j_0\dotso j_\ell} h^{i_k} \dotsm h^{i_0}
    h^{j_\ell} \dotsm h^{j_0} .
  \end{equation}
  We have
  \begin{equation*}
    \om_{i_0\dotso i_k} \om_{j_0\dotso j_\ell} h^{i_k} \dotso h^{i_0}
    h^{j_\ell} \dotsm h^{j_0}
    = (-1)^{k\ell+(k+1)(\ell+1)} \om_{j_0\dotso j_\ell} \om_{i_0\dots
      i_k} h^{j_\ell} \dotsm h^{j_0} h^{i_k} \dotsm h^{i_0} .    
  \end{equation*}
  The expression \eqref{s2} changes sign on exchange of
  $(i_0,\dotsc,i_k)$ and $(j_0,\dotsc,j_\ell)$, and thus vanishes.
  \Endproof

\section{The Maurer-Cartan set of an \Linf-algebra}

\Linf-algebras are a generalization of dg Lie algebras in which the
Jacobi rule is only satisfied up to a hierarchy of higher homotopies.
In this section, we start by recalling the definition of
\Linf-algebras. Following Sullivan \cite{Sullivan} and Hinich
\cite{Hinich}, we represent the homotopy type of an \Linf-algebra $\g$
by the simplicial set $\MC_\bull(\g)=\MC(\g\o\Om_\bull)$. We prove
that this is a Kan complex, and that under certain additional
hypotheses, it is a homotopy invariant of the \Linf-algebra $\g$.

An operation $[x_1,\dotsc,x_k]$ on a graded vector space $\g$ is called
graded antisymmetric if
\begin{equation*} \label{symmetry}
  [x_1,\dotsc,x_i,x_{i+1},\dotsc,x_k] + (-1)^{|x_i||x_{i+1}|}
  [x_1,\dotsc,x_{i+1},x_i,\dotsc,x_k] = 0
\end{equation*}
for all $1\le i\le k-1$. Equivalently, $[x_1,\dotsc,x_k]$ is a linear
map from $\Wedge^k\g$ to $\g$, where $\Wedge^k\g$ is the $k$th
exterior power of the graded vector space $\g$, that is, the $k$th
symmetric power of $s^{-1}\g$.
\begin{definition}
  An \Linf-algebra is a graded vector space $\g$ with a sequence
  $[x_1,\dotsc,x_k]$, $k>0$ of graded antisymmetric operations of
  degree $2-k$, or equivalently, homogeneous linear maps
  $\Wedge^k\g\to\g$ of degree $2$, such that for each $n>0$, the
  $n$-Jacobi rule holds:
  \begin{equation*}
    \sum_{k=1}^n (-1)^k \sum_{\substack{ i_1<\dotsb<i_k ;
    j_1<\dotsb<j_{n-k} \\
    \{i_1,\dotsc,i_k\}\cup\{j_1,\dotsc,j_{n-k}\}=\{1,\dotsc,n\} }}
    (-1)^\eps \, [[x_{i_1},\dotsc,x_{i_k}],x_{j_1},\dotsc,x_{j_{n-k}}] = 0 .
  \end{equation*}
  Here, the sign $(-1)^\eps$ equals the product of the sign $(-1)^\pi$
  associated to the permutation
  \begin{equation*}
    \pi = \bigl( \begin{smallmatrix} 1 & \dotso & k & k+1 & \dotso & n
      \\ i_1 & \dotso & i_k & j_1 & \dotso & j_{n-k}
    \end{smallmatrix} \bigr)
  \end{equation*}
  with the sign associated by the Koszul sign convention to the action
  of $\pi$ on the elements $(x_1,\dotsc,x_n)$ of $\g$.
\end{definition}

In terms of the graded symmetric operations
\begin{equation*}
  \ell_k(y_1,\dotsc,y_k) = (-1)^{\sum_{i=1}^k(k-i+1)|y_i|} \,
  s^{-1}[sy_1,\dotsc,sy_k]
\end{equation*}
of degree $1$ on the graded vector space $s^{-1}\g$, the Jacobi rule
simplifies to become
\begin{equation*}
  \sum_{k=1}^n \sum_{\substack{ i_1<\dotsb<i_k , j_1<\dotsb<j_{n-k} \\
      \{i_1,\dotsc,i_k\}\cup\{j_1,\dotsc,j_{n-k}\}=\{1,\dotsc,n\} }}
  (-1)^{\tilde{\eps}} \,
  \{\{y_{i_1},\dotsc,y_{i_k}\},y_{j_1},\dotsc,y_{j_{n-k}}\} = 0 ,
\end{equation*}
where $(-1)^{\tilde{\eps}}$ is the sign associated by the Koszul sign
convention to the action of $\pi$ on the elements $(y_1,\dotsc,y_n)$ of
$s^{-1}\g$.  This is a small modification of the conventions of Lada
and Markl \cite{LM}: their operations $l_k$ are related to ours by a
sign
\begin{equation*}
  l_k(x_1,\dotsc,x_k) = (-1)^{\binom{k+1}{2}} \, [x_1,\dotsc,x_k] .
\end{equation*}

The operation $x\mapsto[x]$ makes the graded vector space $\g$ into a
cochain complex, by the 1-Jacobi rule $[[x]]=0$. Because of the
special role played by the operation $[x]$, we denote it by $\delta$.
An \Linf-algebra with $[x_1,\dotsc,x_k]=0$ for $k>2$ is the same thing
as a dg Lie algebra. A quasi-isomorphism of \Linf-algebras is a
quasi-isomorphism of the underlying cochain complexes.

The lower central filtration on an \Linf-algebra $\g$ is the canonical
decreasing filtration defined inductively by $F^1\g=\g$ and, for
$i>1$,
\begin{equation*}
  F^i\g = \sum_{i_1+\dotsb+i_k=i} [F^{i_1}\g,\dotsc,F^{i_k}\g] .
\end{equation*}
\begin{definition}
  An \Linf-algebra $\g$ is nilpotent if the lower central series
  terminates, that is, if $F^i\g=0$ for $i\gg0$.
\end{definition}

If $\g$ is a nilpotent \Linf-algebra, the curvature
\begin{equation*}
  \CF(\alpha) = \delta\alpha + \sum_{\ell=2}^\infty \frac{1}{\ell!} \,
  [\alpha^{\wedge \ell}] \in \g^2
\end{equation*}
is defined, and polynomial in $\alpha$. If $\g$ is a dg Lie algebra,
the curvature equals
\begin{equation*}
  \CF(\alpha) = \delta \alpha + \half[\alpha,\alpha] ;
\end{equation*}
this expression is familiar from the theory of connections on
principal bundles.

\begin{definition}
  The Maurer-Cartan set $\MC(\g)$ of a nilpotent \Linf-algebra $\g$ is
  the set of those $\alpha\in\g^1$ satisfying the Maurer-Cartan
  equation
  \begin{equation}
    \label{MaurerCartan}
    \CF(\alpha) = 0 .
  \end{equation}
\end{definition}
An \Linf-algebra is abelian if the bracket $[x_1,\dotsc,x_k]$ vanishes
for $k>1$. In this case, the Maurer-Cartan set is the set of
$1$-cocycles $Z^1(\g)$ of $\g$.

Let $\g$ be a nilpotent \Linf-algebra. For any element
$\alpha\in\g^1$, the formula
\begin{equation*}
  [x_1,\dotsc,x_k]_\alpha = \sum_{\ell=0}^\infty \frac{1}{\ell!} \,
  [\alpha^{\wedge\ell},x_1,\dotsc,x_k]
\end{equation*}
defines a new sequence of brackets on $\g$, where
$[\alpha^{\wedge\ell},x_1,\dotsc,x_k]$ is an abbreviation for
\begin{equation*}
    [\underbrace{\alpha,\dotsc,\alpha}_{\text{$\ell$ times}},x_1,\dotsc,x_k] .
\end{equation*}

\begin{proposition}
  If $\alpha\in\MC(\g)$, then the brackets $[x_1,\dotsc,x_k]_\alpha$
  make $\g$ into an \Linf-algebra.
\end{proposition}
\Proof
  Applying the $(m+n)$-Jacobi relation to the sequence $(\alpha^{\wedge
    m},x_1,\dotsc,x_n)$ and summing over $m$, we obtain the $n$-Jacobi
  relation for the brackets $[x_1,\dotsc,x_k]_\alpha$.
\Endproof

\begin{lemma}
  The curvature satisfies the Bianchi identity
  \begin{equation}
    \label{Bianchi}
    \delta\CF(\alpha) + \sum_{\ell=1}^\infty \frac{1}{\ell!} \,
    [\alpha^{\wedge\ell},\CF(\alpha)] = 0 .
  \end{equation}
\end{lemma}
\Proof
  The $n$-Jacobi relation for $(\alpha^{\wedge n})$ shows that
  \begin{equation*}
    \sum_{\ell=0}^n \frac{1}{\ell!(n-\ell)!} \, [\alpha^{\wedge
      \ell},[\alpha^{\wedge(n-\ell)}]] = 0 .
  \end{equation*}
  Summing over $n>0$, we obtain the lemma.
\Endproof

If $\g$ is an \Linf-algebra and $\Om$ is a dg commutative algebra,
then the tensor product $\g\o\Om$ is an \Linf-algebra, with brackets
\begin{equation*}
  \begin{cases}
    [x\o a] = [x] \o a + (-1)^{|x|} x \o da , & \\[5pt]
    [x_1\o a_1,\dotsc,x_k\o a_k] = (-1)^{\sum_{i<j} |x_i|\,|a_j|}
    [x_1,\dotsc,x_k] \o a_1\dotsm a_k , & k\ne1 .
  \end{cases}
\end{equation*}
The functor $\MC(\g)$ extends to a covariant functor $\MC(\g,\Om) =
\MC(\g\o\Om)$ from dg commutative algebras to sets, that is, a
presheaf on the category of dg affine schemes over $\K$. If $X_\bull$
is a simplicial set, we have
\begin{equation*}
  \MC(\g,\Om(X_\bull)) \cong \ss(X_\bull,\MC_\bull(\g)) .
\end{equation*}

If $\g$ is a nilpotent \Linf-algebra, let $\MC_\bull(\g)$ be the
simplicial set
\begin{equation*}
  \MC_\bull(\g) = \MC(\g,\Om_\bull) .
\end{equation*}
In other words, the $n$-simplices of $\MC_\bull(\g)$ are differential
forms $\alpha$ on the $n$-simplex $\DELTA^n$, of the form
\begin{equation*}
  \alpha = \sum_{i=0}^n \alpha_i
\end{equation*}
where $\alpha_i\in\g^{1-i}\o\Om^i(\DELTA^n)$, such that
\begin{equation} \label{mc}
    (d+\delta)\alpha + \sum_{\ell=2}^\infty \frac{1}{\ell!} \,
    [\alpha^{\wedge\ell}] = 0 .    
\end{equation}
Before developing the properties of this functor, we recall how it
emerges naturally from Sullivan's approach \cite{Sullivan} to rational
homotopy theory.

If $\g$ is an \Linf-algebra which is finite-dimensional in each degree
and bound\-ed below, we may associate to it the dg commutative algebra
$C^*(\g)$ of cochains. The underlying graded commutative algebra of
$C^*(\g)$ is $\Wedge\g^\vee=S(\g[1]^\vee)$, the free graded
commutative algebra on the graded vector space $\g[1]^\vee$ which
equals $\bigl(\g^{1-i}\bigr){}^\vee$ in degree $i$.  The differential
$\delta$ of $C^*(\g)$ is determined by its restriction to the space of
generators $\g[1]^\vee\subset C^*(\g)$, on which it equals the sum
over $k$ of the adjoints of the operations $\ell_k$. The resulting
graded derivation satisfies the equation $\delta^2=0$ if and only if
$\g$ is an \Linf-algebra.

As explained in the introduction, the simplicial set
\begin{equation*}
  \Spec_\bull(\CA) = \dga(\CA,\Om_\bull) .
\end{equation*}
may be viewed as an analogue in homotopical algebra of the spectrum of
a commutative algebra. Applied to $C^*(\g)$, we obtain a simplicial
set $\Spec_\bull(C^*(\g))$ which has a natural identification with the
simplicial set $\MC_\bull(\g)$.

The homotopy groups of a nilpotent \Linf-algebra $\g$ are defined as
\begin{equation*}
  \pi_i(\g) = \pi_i(\MC_\bull(\g)) .
\end{equation*}
In particular, the set of components $\pi_0(\g)$ of $\g$ is the
quotient of $\MC(\g)$ by the nilpotent group associated to the
nilpotent Lie algebra $\g^0$. This plays a prominent role in
deformation theory: it is the moduli set of deformations of $\g$.

In order to establish that $\MC_\bull(\g)$ is a Kan complex, we use
the Poincar\'e lemma. Let $0\le i\le n$. By \eqref{Poincare}, we see
that
\begin{equation*}
  \Id_n = \eps^i_n + (d+\delta)h^i_n + h^i_n(d+\delta) .
\end{equation*}
If $\alpha\in\MC_n(\g)$, we see that
\begin{align*}
  \alpha &= \eps^i_n\alpha + (d+\delta)h^i_n\alpha +
  h^i_n(d+\delta)\alpha \\
  &= \eps^i_n\alpha + R^i_n\alpha - \sum_{\ell=2}^\infty
  \frac{1}{\ell!} \, h^i_n[\alpha^{\wedge\ell}] ,
\end{align*}
where $R^i_n=(d+\delta)h_n^i$. Introduce the space
\begin{equation*}
  \mc_n(\g) = \{ (d+\delta)\alpha \mid \alpha \in \bigl(\g\o\Om
  \bigr){}^0 \bigr\} .
\end{equation*}
\begin{lemma} \label{Unique}
  Let $\g$ be a nilpotent \Linf-algebra. The map
  $\alpha\mapsto(\eps^i_n\alpha,R^i_n\alpha)$ induces an isomorphism
  between $\MC_n(\g)$ and $\MC(\g)\times\mc_n(\g)$.
\end{lemma}
\Proof
  Given $\mu\in\MC(\g)$ and $\nu\in\mc_n(\g)$, let $\alpha_0=\mu+\nu$
  and define differential forms $(\alpha_k)_{k>0}$ inductively by the
  formula
  \begin{equation} \label{iterateMC}
    \alpha_{k+1} = \alpha_0 - \sum_{\ell=2}^\infty \frac{1}{\ell!} \,
    h^i_n[\alpha_k^{\wedge\ell}] .
  \end{equation}
  Then for all $k$, we have $\eps^i_n\alpha_k=\mu$ and
  $R^i_n\alpha_k=\nu$. The sequence is eventually constant, since by
  induction, we see that
  \begin{align*}
    \alpha_{k+1} - \alpha_k &= \sum_{\ell=2}^\infty \frac{1}{\ell!} \,
    \sum_{j=1}^\ell h^i_n \bigl[ \alpha_{k-1}^{\wedge
      j-1},\alpha_{k-1}-\alpha_k,\alpha_k^{\wedge\ell-j} \bigr] \\
    &\in F^{k+1}\g \o \Om_n .
  \end{align*}

  The limit
  \begin{equation*}
    \alpha = \lim_{k\to\infty} \alpha_k
  \end{equation*}
  satisfies
  \begin{equation*}
    \alpha = \alpha_0 - \sum_{\ell=2}^\infty \frac{1}{\ell!} \,
    h^i_n[\alpha^{\wedge\ell}] .
  \end{equation*}
  Applying the operator $d+\delta$, we see that
  \begin{equation*}
    (d+\delta)\alpha = \delta\mu - \sum_{\ell=2}^\infty
    \frac{1}{\ell!} \, (d+\delta)h^i_n[\alpha^{\wedge\ell}] ,
  \end{equation*}
  and hence that
  \begin{align*}
    \CF(\alpha) &= \delta\mu + \sum_{\ell=2}^\infty
    \frac{1}{\ell!} \, [\alpha^{\wedge\ell}] - \sum_{\ell=2}^\infty
    \frac{1}{\ell!} \, (d+\delta)h^i_n[\alpha^{\wedge\ell}] \\
    &= \CF(\mu) + \sum_{\ell=2}^\infty \frac{1}{\ell!} \,
    h^i_n(d+\delta)[\alpha^{\wedge\ell}] \\
    &= \CF(\mu) + h^i_n(d+\delta)\CF(\alpha) .
  \end{align*}
  The Bianchi identity \eqref{Bianchi} implies that
  \begin{align*}
    \CF(\alpha) &= \CF(\mu) - \sum_{\ell=1}^\infty \frac{1}{\ell!}
    h^i_n[\alpha^{\wedge\ell},\CF(\alpha)] \\
    &= \sum_{\ell=1}^\infty \frac{1}{\ell!}
    h^i_n[\alpha^{\wedge\ell},\CF(\alpha)] .
  \end{align*}
  The nilpotence of $\g$ implies that $\CF(\alpha)=0$; it follows that
  $\alpha$ is an element of $\MC_n(\g)$ with $\eps^i_n\alpha=\mu$ and
  $R^i_n\alpha=\nu$.
  
  If $\alpha$ and $\beta$ are a pair of elements of $\MC_n(\g)$ such
  that $\eps^i_n\alpha=\eps^i_n\beta$ and $R^i_n\alpha=R^i_n\beta$,
  then
  \begin{equation*}
    \alpha - \beta = - \sum_{\ell=2}^\infty \frac{1}{\ell!} \,
    \sum_{j=1}^\ell h^i_n \bigl[ \alpha^{\wedge
      j-1},\alpha-\beta,\beta^{\wedge\ell-j} \bigr] .
  \end{equation*}
  This shows, by induction, that $\alpha-\beta\in F^i\g$ for all
  $i>0$, and hence, by the nilpotence of $\g$, that $\alpha=\beta$.
\Endproof

The following result is due to Hinich \cite{Hinich} when $\g$ is a dg
Lie algebra.
\begin{proposition} \label{fibration}
  If $f:\g\to\h$ is a surjective morphism of nilpotent \Linf-algebras,
  the induced morphism
  \begin{equation*}
    \MC_\bull(f):\MC_\bull(\g)\to\MC_\bull(\h)
  \end{equation*}
  is a fibration of simplicial sets.
\end{proposition}
\Proof
  Let $0\le i\le n$. Given a horn
  \begin{equation*} 
    \beta \in \ss(\Lambda^n_i,\MC_\bull(\g))
  \end{equation*}
  and an $n$-simplex $\gamma\in\MC_n(\h)$ such that
  $\p_j\gamma=f(\p_j\beta)$ for $j\ne i$, we wish to construct an
  element $\alpha\in f^{-1}(\gamma)\subset\MC_n(\g)$ such that
  \begin{equation*}
    \p_j\alpha=\p_j\beta
  \end{equation*}
  for $j\ne i$.
  
  Since $f:\g\o\Om_\bull\to\h\o\Om_\bull$ is a Kan fibration, there
  exists an extension $\rho\in\g\o\Om_n$ of $\beta$ of total degree
  $1$ such that $f(\rho)=\alpha$. Let $\alpha$ be the unique element
  of $\MC_n(\g)$ such that $\eps^i_n\alpha=\eps^i_n\rho$ and
  $R^i_n\alpha=R^i_n\rho$. If $j\ne i$, we have
  $\eps^i_n\p_j\alpha=\eps^i_n\p_j\beta$ and
  $R^i_n\p_j\alpha=R^i_n\p_j\beta$ and hence, by Lemma \ref{Unique},
  $\p_j\alpha=\p_j\beta$; thus, $\alpha$ fills the horn $\beta$. We
  also have $f(\eps^i_n\alpha)=f(\eps^i_n\rho)=\eps^i_n\gamma$ and
  $f(R^i_n\alpha)=f(R^i_n\rho)=R^i_n\gamma$, hence $f(\alpha)=\gamma$.
\Endproof

The category of nilpotent \Linf-algebras concentrated in degrees
$(-\infty,0]$ is a variant of Quillen's model for rational homotopy of
nilpotent spaces \cite{Quillen}. By the following theorem, the functor
$\MC_\bull(\g)$ carries quasi-isomorphisms of such \Linf-algebras to
homotopy equivalences of simplicial sets.
\begin{theorem}
  If $\g$ and $\h$ are \Linf-algebras concentrated in degrees
  $(-\infty,0]$, and $f:\g\to\h$ is a quasi-isomorphism, then
  \begin{equation*}
    \MC_\bull(f) : \MC_\bull(\g) \to \MC_\bull(\h)
  \end{equation*}
  is a homotopy equivalence.
\end{theorem}
\Proof
  Filter $\g$ by \Linf-algebras $F^j\g$, where
  \begin{align*}
    (F^{2j}\g)^i &=
    \begin{cases}
      0 , & i+j>0 , \\
      Z^{-j}(\g) , & i+j=0 , \\
      \g^i , & i+j<0 ,
    \end{cases} &
    (F^{2j+1}\g)^i &=
    \begin{cases}
      0 , & i+j>0 , \\
      B^{-j}(\g) , & i+j=0 , \\
      \g^i , & i+j<0 ,
    \end{cases}
  \end{align*}
  and similarly for $\h$. If $j>k$, there is a morphism of fibrations
  of simplicial sets
  \begin{equation*}
    \begin{CD}
      \MC_\bull(F^j\g) @>>> \MC_\bull(F^k\g) @>>>
      \MC_\bull(F^k\g/F^j\g) \\
      @VVV @VVV @VVV \\
      \MC_\bull(F^j\h) @>>> \MC_\bull(F^k\h) @>>>
      \MC_\bull(F^k\h/F^j\h)
    \end{CD}
  \end{equation*}
  We have
  \begin{equation*}
    \MC_\bull(F^{2j}\g/F^{2j+1}\g) \cong  \MC_\bull(H^{-j}(\g)) \cong
    \MC_\bull(H^{-j}(\h)) \cong \MC_\bull(F^{2j}\h/F^{2j+1}\h) .
  \end{equation*}
  The simplicial sets
  \begin{align*}
    \MC_\bull(F^{2j+1}\g/F^{2j+2}\g) &\cong B^{-j}(\g) \o
    \Om_\bull^{j+1} , \quad \text{and} \\
    \MC_\bull(F^{2j+1}\h/F^{2j+2}\h) &\cong B^{-j}(\h) \o
    \Om_\bull^{j+1}
  \end{align*}
  are contractible by Lemma \ref{contractible}. The proposition follows.
\Endproof

Let $\m$ be a nilpotent commutative ring; that is, $\m^{\ell+1}=0$ for
some $\ell$. If $\g$ is an \Linf-algebra, then $\g\o\m$ is nilpotent;
this is the setting of of formal deformation theory. In this context
too, the functor $\MC_\bull(\g,\m)=\MC_\bull(\g\o\m)$ takes
quasi-isomorphisms of \Linf-algebras to homotopy equivalences of
simplicial sets.

\begin{proposition}
  If $f:\g\to\h$ is a quasi-isomorphism of \Linf-algebras and $\m$ is
  a nilpotent commutative ring, then
  \begin{equation*}
    \MC_\bull(f,\m) : \MC_\bull(\g,\m) \to \MC_\bull(\h,\m) 
  \end{equation*}
  is a homotopy equivalence.
\end{proposition}
\Proof
  We argue by induction on the nilpotence length $\ell$ of $\m$. There
  is a morphism of fibrations of simplicial sets
  \begin{equation*}
    \begin{CD}
      \MC_\bull(\g,\m^2) @>>> \MC_\bull(\g,\m) @>>>
      \MC_\bull(\g\o\m/\m^2) \\
      @VVV @VVV @VVV \\
      \MC_\bull(\h,\m^2) @>>> \MC_\bull(\h,\m) @>>>
      \MC_\bull(\h\o\m/\m^2)
    \end{CD}
  \end{equation*}
  The abelian \Linf-algebras $\g\o\m/\m^2$ and $\h\o\m/\m^2$ are
  quasi-isomorphic, hence the morphism
  \begin{equation*}
    \MC_\bull(\g\o\m/\m^2) \to \MC_\bull(\h\o\m/\m^2)
  \end{equation*}
  is a homotopy equivalence. The result follows by induction on
  $\ell$.
\Endproof

\section{The functor $\GAMMA_\bull(\g)$}

In this section, we study the functor $\GAMMA_\bull(\g)$; we prove
that it is homotopy equivalent to $\MC_\bull(\g)$, and show that it
specializes to the Deligne groupoid when $\g$ is concentrated in
degrees $[0,\infty)$. Fix a gauge $s_\bull$, for example Dupont's
operator \eqref{s}.

The simplicial set $\GAMMA_\bull(\g)$ associated to a nilpotent
\Linf-algebra is the simplicial subset of $\MC_\bull(\g)$ consisting
of those Maurer-Cartan forms annihilated by $s_\bull$:
\begin{equation} \label{gamma}
  \GAMMA_\bull(\g) = \{ \alpha \in \MC_\bull(\g) \mid s_\bull\alpha =
  0 \} .
\end{equation}
For any simplicial set $X_\bull$, the set of simplicial maps
$\ss(X_\bull,\GAMMA_\bull(\g))$ equals the set of Maurer-Cartan
elements $\alpha\in\MC(\g,X_\bull)$ such that $s_\bull\alpha=0$. This
is reminiscent of gauge conditions, such as the Coulomb gauge, in
gauge theory.

\begin{proposition} \label{DoldKan}
  For abelian $\g$, there is a natural isomorphism
  $\GAMMA_\bull(\g)\cong K_\bull(\g[1])$.
\end{proposition}
\Proof
  If $\alpha\in\GAMMA_n(\g)$, then $(d+\delta)\alpha=s_n\alpha=0$,
  hence by \eqref{Dupont},
  \begin{equation*}
    \alpha = P_n\alpha + s_n(d+\delta)\alpha + (d+\delta)s_n\alpha =
    P_n\alpha .
  \end{equation*}
  Thus $\GAMMA_n(\g)\subset K_n(\g[1])$. Conversely, if $\alpha\in
  K_n(\g[1])$, then $P_n\alpha=\alpha$, hence $s_n\alpha=0$. Thus
  $K_n(\g[1])\subset\GAMMA_n(\g)$.
\Endproof

We show that $\GAMMA_\bull(\g)$ is an $\infty$-groupoid, and in
particular, a Kan complex: the heart of the proof is an iteration,
similar to the iteration \eqref{iterateMC}, which solves the
Maurer-Cartan equation on the $n$-simplex $\Delta^n$ in the gauge
$s_n\alpha=0$.

\begin{definition}
  An $n$-simplex $\alpha\in\GAMMA_n(\g)$ is \bold{thin} if
  $I_{0\dotso n}(\alpha)=0$.
\end{definition}

\begin{lemma} \label{unique}
  Let $\g$ be a nilpotent \Linf-algebra. The map
  $\alpha\mapsto(\eps^i_n\alpha,P_nR^i_n\alpha)$ induces an
  isomorphism between $\GAMMA_n(\g)$ and
  $\MC(\g)\times P_n[\mc_n(\g)]$.
\end{lemma}
\Proof
  Let $0\le i\le n$. By \eqref{Dupont}, we see that
  \begin{align*}
    \Id_n &= P_n + (d+\delta)s_n + s_n(d+\delta) \\
    &= \eps^i_n + (d+\delta)(P_nh^i_n+s_n) + (P_nh^i_n+s_n)(d+\delta)
    .
  \end{align*}
  It follows that if $\alpha\in\GAMMA_n(\g)$,
  \begin{equation} \label{iterate}
    \alpha = \eps^i_n\alpha + P_nR^i_n\alpha - \sum_{\ell=2}^\infty
    \frac{1}{\ell!} \, (P_nh^i_n+s_n)[\alpha^{\wedge\ell}] .
  \end{equation}
  
  Given $\mu\in\MC(\g)$ and $\nu\in P_n[\mc_n(\g)]$, let
  $\alpha_0=\mu+\nu$ and define differential forms $(\alpha_k)_{k>0}$
  inductively by the formula
  \begin{equation*}
    \alpha_k = \alpha_0 - \sum_{\ell=2}^\infty \frac{1}{\ell!} \,
    (P_nh^i_n+s_n)[\alpha_{k-1}^{\wedge\ell}] .
  \end{equation*}
  Then for all $k$, we have $s_n\alpha_k=0$, $\eps^i_n\alpha_k=\mu$
  and $P_nR^i_n\alpha_k=\nu$. The sequence $(\alpha_k)$ is eventually
  constant, since by induction, we see that
  \begin{align*}
    \alpha_k - \alpha_{k-1} &= \sum_{\ell=2}^\infty \frac{1}{\ell!} \,
    \sum_{j=1}^\ell (P_nh^i_n+s_n) \bigl[ \alpha_{k-2}^{\wedge j-1} ,
    \alpha_{k-2}-\alpha_{k-1},\alpha_{k-1}^{\wedge\ell-j} \bigr] \\
    &\in F^k\g \o \Om_n .
  \end{align*}

  The limit
  \begin{equation*}
    \alpha = \lim_{k\to\infty} \alpha_k
  \end{equation*}
  satisfies
  \begin{equation*}
    \alpha = \alpha_0 - \sum_{\ell=2}^\infty \frac{1}{\ell!} \,
    (Ph^i_n+s_n)[\alpha^{\wedge\ell}] .
  \end{equation*}
  By the same argument as in the proof of Lemma~\ref{Unique}, it
  follows that
  \begin{align*}
    \CF(\alpha) &= \CF(\mu) - \sum_{\ell=1}^\infty \frac{1}{\ell!}
    (Ph^i_n+s_n)[\alpha^{\wedge\ell},\CF(\alpha)] \\
    &= \sum_{\ell=1}^\infty \frac{1}{\ell!}
    (Ph^i_n+s_n)[\alpha^{\wedge\ell},\CF(\alpha)] .
  \end{align*}
  The nilpotence of $\g$ implies that $\CF(\alpha)=0$; it follows that
  $\alpha$ is an element of $\gamma_n(\g)$ with $\eps^i_n\alpha=\mu$
  and $PR^i_n\alpha=\nu$.
  
  If $\alpha$ and $\beta$ are a pair of elements of $\GAMMA_n(\g)$
  such that $\eps^i_n\alpha=\eps^i_n\beta$ and
  $P_nR^i_n\alpha=P_nR^i_n\beta$, then
  \begin{equation*}
    \alpha - \beta = - \sum_{\ell=2}^\infty \frac{1}{\ell!} \,
    \sum_{j=1}^\ell (P_nh^i_n+s_n) \bigl[ \alpha^{\wedge
      j-1},\alpha-\beta,\beta^{\wedge\ell-j} \bigr] .
  \end{equation*}
  This shows, by induction, that $\alpha-\beta\in F^i\g$ for all
  $i>0$, and hence, by the nilpotence of $\g$, that $\alpha=\beta$.
\Endproof

\begin{theorem} \label{GAMMA}
  If $\g$ is a nilpotent \Linf-algebra, $\GAMMA_\bull(\g)$ is an
  $\infty$-groupoid. If $\g$ is concentrated in degrees
  $(-\ell,\infty)$, respectively $(-\ell,0]$, then $\GAMMA_\bull(\g)$
  is an $\ell$-groupoid, resp.\ an $\ell$-group.
\end{theorem}
\Proof
  Let
  \begin{equation*} 
    \beta \in \ss(\Lambda^n_i,\GAMMA_\bull(\g))
  \end{equation*}
  be a horn in $\GAMMA_\bull(\g)$. The differential form
  \begin{equation*}
    \alpha_0 = \eps^i_n\beta + (d+\delta) \sum_{k=1}^{n-1}
    \sum_{\substack{i_1<\dots<i_k \\ i\notin\{i_1,\dotsc,i_k\} }}
    \om_{i_1\dotso i_k} \o I_{ii_1\dotso i_k}(\beta) \in \MC(\g) \times
    P_n[\mc_n(\g)]
  \end{equation*}
  satisfies $I_{0\dotso n}(\alpha_0)=0$. The solution
  $\alpha\in\GAMMA_n(\g)$ of the equation
  \begin{equation*}
    \alpha = \alpha_0 - \sum_{\ell=2}^\infty \frac{1}{\ell!} \,
    (P_nh^i_n+s_n)[\alpha^{\wedge\ell}]
  \end{equation*}
  constructed in Lemma \ref{unique} is thin and
  $\xi^n_i(\alpha)=\beta$. Thus $\GAMMA_\bull(\g)$ is an
  $\infty$-groupoid.
  
  If $\g^{1-n}=0$, it is clear that every $n$-simplex
  $\alpha\in\GAMMA_n(\g)$ is thin, while if $\g^1=0$, then
  $\GAMMA_\bull(\g)$ is reduced.
\Endproof

Given $\mu\in\MC(\g)$ and $x_{i_1\dotso i_k}\in\g^{1-k}$, $1\le
i_1<\dots<i_k\le n$, let $\alpha_n^\mu(x_{i_1\dotso i_k}) \in
\GAMMA_n(\g)$ be the solution of \eqref{iterate} with
$\eps_n^0\alpha_n^\mu(x_{i_1\dotso i_k})=\mu$ and
\begin{equation*}
  R_n^0\alpha_n^\mu(x_{i_1\dotso i_k}) = \sum_{k=1}^n \sum_{1\le
    i_1<\dots<i_k\le n} \om_{i_1\dotso i_k} \o x_{i_1\dotso i_k} .
\end{equation*}

\begin{definition}
  The $n$th generalized Campbell-Hausdorff series associated to the
  gauge $s_\bull$ is the function of $\mu\in\MC(\g)$ and
  \begin{equation*}
    x_{i_1\dots i_k}\in\g^{1-k} , \quad 1\le i_1<\dots<i_k\le n ,    
  \end{equation*}
  given by the formula
  \begin{equation*}
    \rho_n^\mu(x_{i_1\dotso i_k}) = I_{1\dotso n} \bigl(
    \alpha_n^\mu(x_{i_1\dotso i_k}) \bigr) \in \g^{2-n} .
  \end{equation*}
\end{definition}

If $\g$ is concentrated in degrees $(-\infty,0]$, then the
Maurer-Cartan element $\mu$ equals $0$, and is omitted from the
notation for $\alpha_n(x_{i_1\dotso i_k})$ and $\rho_n(x_{i_1\dots
  i_k})$.

Since $\alpha_2^\mu(x_1,x_2,x_{12})$ is a flat connection $1$-form on
the $2$-simplex, its monodromy around the boundary must be trivial.
(The $2$-simplex is simply connected.) In terms of the generalized
Campbell-Hausdorff series $\rho_2^\mu(x_1,x_2,x_{12})$, this gives the
equation
\begin{equation*}
  e^{x_1} = e^{\rho_2^\mu(x_1,x_2,x_{12})} e^{x_2}
\end{equation*}
in the Lie group associated to the nilpotent Lie algebra $\g^0$. Thus,
the simplicial set $\GAMMA_\bull(\g)$ (indeed, its $2$-skeleton)
determines $\rho_2^\mu(x_1,x_2,x_{12})$ as a function of $x_1$, $x_2$
and $x_{12}$. In the Dupont gauge, modulo terms involving more than
two brackets, it equals
\begin{align*}
  \rho_2^\mu(x_1,x_2,x_{12}) &= x_1 - x_2 + \half [x_1,x_2]_\mu +
  \half [x_{12}]_\mu \\
  &+ \tfrac{1}{12} [x_1+x_2,[x_1,x_2]_\mu]_\mu +
  \tfrac{1}{6} [[x_1+x_2]_\mu,x_1,x_2]_\mu \\
  &+ \tfrac{1}{6} [[x_1+x_2]_\mu,x_{12}]_\mu - \tfrac{1}{12}
  [x_1+x_2,[x_{12}]_\mu]_\mu + \dotsb .
\end{align*}

If $\g$ is a dg Lie algebra, the thin $2$-simplices define a
composition on the $1$-simplices of $\GAMMA_\bull(\g)$ which is
strictly associative.
\begin{proposition}
  If $\g$ is a dg Lie algebra, the composition
  \begin{equation*}
    \rho_2^\mu(x_1,x_2) : \g \o \g \to \g
  \end{equation*}
  is associative.
\end{proposition}
\Proof
  It suffices to show that $\rho_3^\mu(x_1,x_2,x_3,x_{ij}=0)=0$, in
  other words, if three faces of a thin 3-simplex are thin, then the
  fourth is. The iteration leading to the solution $\alpha$ of
  \eqref{iterate} with initial conditions
  \begin{equation*}
    \alpha_0 = \mu + (d+\delta) (t_1x_1 + t_2x_2 + t_3x_3)
  \end{equation*}
  lies in the space $\g^0\o\Om^1_3\oplus\g^1\o\Om^0_3$, hence
  $I_{123}(\alpha)=0$.
\Endproof

In particular, if $\g$ is a dg Lie algebra concentrated in degrees
$(-2,\infty)$, $\GAMMA_\bull(\g)$ is the nerve of a strict
$2$-groupoid, that is, a groupoid enriched in groupoids; in this way,
we see that $\GAMMA_\bull(\g)$ generalizes the Deligne $2$-groupoid
(Deligne \cite{Deligne}, Getzler \cite{darboux}).

Although it is not hard to derive explicit formulas for the
generalized Campbell-Hausdorff series up to any order, we do not know
any closed formulas for them except when $n=1$, in which case it is
independent of the gauge. We now derive a closed formula for
$\rho_1^\mu(x)$, which resembles Cayley's famous formula for the
series solution of the ordinary differential equation $x'(t)=f(x(t))$.

To each rooted tree, associate the word obtained by associating to a
vertex with $i$ branches the operation $[x,a_1,\dotsc,a_i]_\mu$.
Multiply the resulting word by the number of total orders on the
vertices of the tree such that each vertex precedes its parent. Let
$\e_\mu^k(x)$ be the sum of these terms over all rooted trees with $k$
vertices. For example, $\e^1_\mu(x)=[x]_\mu$,
$\e^2_\mu(x)=[x,[x]_\mu]_\mu$ and
\begin{equation*}
  \e^3_\mu(x) = [x,[x,[x]_\mu]_\mu]_\mu+[x,[x]_\mu,[x]_\mu]_\mu .
\end{equation*}
The coefficient of a tree $T$ in $\e^k_\mu(x)$ equals the number of
monotone orderings of its vertices, that is, total orderings such that
each vertex is greater than its parent. See Figure~2 for the trees
contributing to $\e^k_\mu(X)$ for small values of $k$.

\begin{figure}
  \begin{align*}
  \e_1(X) &= \begin{picture}(30,40)(-20,12)
    \put(0,15){\circle*{4}}
  \end{picture} \qquad
  \e_2(X) = \begin{picture}(30,40)(-20,12)
    \put(0,10){\circle*{4}}
    \put(0,10){\line( 0,1){ 15}}
    \put(0,25){\circle*{4}}
  \end{picture} \qquad
  \e_3(X) =
  \begin{picture}(30,40)(-20,12)
    \put(0,0){\circle*{4}}
    \put(0,0){\line( 0,1){ 15}}
    \put(0,15){\circle*{4}}
    \put(0,15){\line( 0,1){ 15}}
    \put(0,30){\circle*{4}}
  \end{picture} +
  \begin{picture}(30,40)(-20,12)
    \put(-7.5,15){\circle*{4}}
    \put(-7.5,15){\line( 1,2){ 7.5}}
    \put(7.5,15){\circle*{4}}
    \put(7.5,15){\line( -1,2){ 7.5}}
    \put(0,30){\circle*{4}}
  \end{picture} \\
  \e_4(X) &=
  \begin{picture}(40,60)(-20,25)
    \put(0,0){\circle*{4}}
    \put(0,0){\line( 0,1){ 15}}
    \put(0,15){\circle*{4}}
    \put(0,15){\line( 0,1){ 15}}
    \put(0,30){\circle*{4}}
    \put(0,30){\line( 0,1){ 15}}
    \put(0,45){\circle*{4}}
  \end{picture} +
  \begin{picture}(40,60)(-20,25)
    \put(-7.5,15){\circle*{4}}
    \put(-7.5,15){\line( 1,2){ 7.5}}
    \put(7.5,15){\circle*{4}}
    \put(7.5,15){\line( -1,2){ 7.5}}
    \put(0,30){\circle*{4}}
    \put(0,30){\line( 0,1){ 15}}
    \put(0,45){\circle*{4}}
  \end{picture} + 3
  \begin{picture}(40,60)(-20,25)
    \put(-7.5,15){\circle*{4}}
    \put(-7.5,15){\line( 0,1){ 15}}
    \put(-7.5,30){\circle*{4}}
    \put(-7.5,30){\line( 1,2){ 7.5}}
    \put(7.5,30){\circle*{4}}
    \put(7.5,30){\line( -1,2){ 7.5}}
    \put(0,45){\circle*{4}}
  \end{picture} +
  \begin{picture}(40,60)(-20,25)
    \put(-15,30){\circle*{4}}
    \put(-15,30){\line( 1,1){ 15}}
    \put(15,30){\circle*{4}}
    \put(15,30){\line( -1,1){ 15}}
    \put(0,30){\circle*{4}}
    \put(0,30){\line( 0,1){ 15}}
    \put(0,45){\circle*{4}}
  \end{picture} \\
  \e_5(X) &=
  \begin{picture}(40,80)(-20,40)
    \put(0,0){\circle*{4}}
    \put(0,0){\line( 0,1){ 15}}
    \put(0,15){\circle*{4}}
    \put(0,15){\line( 0,1){ 15}}
    \put(0,30){\circle*{4}}
    \put(0,30){\line( 0,1){ 15}}
    \put(0,45){\circle*{4}}
    \put(0,45){\line( 0,1){ 15}}
    \put(0,60){\circle*{4}}
  \end{picture} +
  \begin{picture}(40,80)(-20,40)
    \put(-7.5,15){\circle*{4}}
    \put(-7.5,15){\line( 1,2){ 7.5}}
    \put(7.5,15){\circle*{4}}
    \put(7.5,15){\line( -1,2){ 7.5}}
    \put(0,30){\circle*{4}}
    \put(0,30){\line( 0,1){ 15}}
    \put(0,45){\circle*{4}}
    \put(0,45){\line( 0,1){ 15}}
    \put(0,60){\circle*{4}}
  \end{picture} + 3
  \begin{picture}(40,80)(-20,40)
    \put(-7.5,15){\circle*{4}}
    \put(-7.5,15){\line( 0,1){ 15}}
    \put(-7.5,30){\circle*{4}}
    \put(-7.5,30){\line( 1,2){ 7.5}}
    \put(7.5,30){\circle*{4}}
    \put(7.5,30){\line( -1,2){ 7.5}}
    \put(0,45){\circle*{4}}
    \put(0,45){\line( 0,1){ 15}}
    \put(0,60){\circle*{4}}
  \end{picture} +
  \begin{picture}(40,80)(-20,40)
    \put(-15,30){\circle*{4}}
    \put(-15,30){\line( 1,1){ 15}}
    \put(15,30){\circle*{4}}
    \put(15,30){\line( -1,1){ 15}}
    \put(0,30){\circle*{4}}
    \put(0,30){\line( 0,1){ 15}}
    \put(0,45){\circle*{4}}
    \put(0,45){\line( 0,1){ 15}}
    \put(0,60){\circle*{4}}
  \end{picture} + 4
  \begin{picture}(40,80)(-20,40)
    \put(-7.5,30){\circle*{4}}
    \put(-7.5,30){\line( 1,2){ 15}}
    \put(7.5,30){\circle*{4}}
    \put(7.5,30){\line( -1,2){ 7.5}}
    \put(0,45){\circle*{4}}
    \put(15,45){\circle*{4}}
    \put(7.5,60){\line( 1,-2){ 7.5}}
    \put(7.5,60){\circle*{4}}
  \end{picture} \\
  &\quad + 4
  \begin{picture}(40,80)(-20,40)
    \put(-7.5,15){\circle*{4}}
    \put(-7.5,15){\line( 0,1){ 15}}
    \put(-7.5,30){\circle*{4}}
    \put(-7.5,30){\line( 0,1){ 15}}
    \put(-7.5,45){\circle*{4}}
    \put(-7.5,45){\line( 1,2){ 7.5}}
    \put(7.5,45){\circle*{4}}
    \put(7.5,45){\line( -1,2){ 7.5}}
    \put(0,60){\circle*{4}}
  \end{picture} + 3
  \begin{picture}(40,80)(-20,40)
    \put(7.5,30){\circle*{4}}
    \put(7.5,30){\line( 0,1){ 15}}
    \put(-7.5,30){\circle*{4}}
    \put(-7.5,30){\line( 0,1){ 15}}
    \put(-7.5,45){\circle*{4}}
    \put(-7.5,45){\line( 1,2){ 7.5}}
    \put(7.5,45){\circle*{4}}
    \put(7.5,45){\line( -1,2){ 7.5}}
    \put(0,60){\circle*{4}}
  \end{picture} + 6
    \begin{picture}(50,80)(-25,40)
    \put(-15,30){\circle*{4}}
    \put(-15,30){\line( 0,1){ 15}}
    \put(-15,45){\circle*{4}}
    \put(-15,45){\line( 1,1){ 15}}
    \put(15,45){\circle*{4}}
    \put(15,45){\line( -1,1){ 15}}
    \put(0,45){\circle*{4}}
    \put(0,45){\line( 0,1){ 15}}
    \put(0,60){\circle*{4}}
  \end{picture} +
  \begin{picture}(60,80)(-30,40)
    \put(23,45){\circle*{4}}
    \put(23,45){\line( -3,2){ 23}}
    \put(-23,45){\circle*{4}}
    \put(-23,45){\line( 3,2){ 23}}
    \put(-7.5,45){\circle*{4}}
    \put(-7.5,45){\line( 1,2){ 7.5}}
    \put(7.5,45){\circle*{4}}
    \put(7.5,45){\line( -1,2){ 7.5}}
    \put(0,60){\circle*{4}}
  \end{picture} \\
\end{align*}
  \caption{Trees representing $\e^k_\mu(X)$ for $k\le5$}
\end{figure}

\begin{proposition}
  \label{exp}
  The $1$-simplex $\alpha_1^\mu(x)\in\GAMMA_1(\g)$ determined by
  $\alpha\in\MC(\g)$ and $x\in\g^0$ is given by the formula
  \begin{equation*}
    \alpha_1^\mu(x) = \alpha - \sum_{k=1}^\infty \frac{t_0^k}{k!} \,
    \e^k_\mu(x) + x\,dt_0 .
  \end{equation*}
\end{proposition}
\Proof
  To show that $\alpha_1^\mu(x)\in\GAMMA_1(\g)$, we must show that it
  satisfies the Maurer-Cartan equation. Let
  \begin{equation*}
    \alpha(t) = \alpha - \sum_{k=1}^\infty \frac{t^k}{k!} \,
    \e_\mu^k(x) .
  \end{equation*}
  It must be shown that
  \begin{equation*}
    \alpha'(t) + \sum_{n=0}^\infty \frac{1}{n!} \, [\alpha(t)^{\wedge
      n},x] = 0 ,
  \end{equation*}
  in other words, that
  \begin{align*}
    \e_\mu^{k+1}(x) &= \sum_{n=0}^\infty \frac{(-1)^n}{n!}
    \sum_{k_1+\dotsb+k_n=k} \frac{k!}{k_1!\dotsm k_n!} \,
    [\e_\mu^{k_1}(x),\dotsc,\e_\mu^{k_n}(x),x]_\alpha \\
    &= \sum_{n=0}^\infty \frac{1}{n!}  \sum_{k_1+\dotsb+k_n=k}
    \frac{k!}{k_1!\dotsm k_n!} \,
    [x,\e_\mu^{k_1}(x),\dotsc,\e_\mu^{k_n}(x)]_\alpha .
  \end{align*}
  This is easily proved by induction on $k$.
\Endproof

Proposition \ref{exp} implies the following formula for the
generalized Campbell-Hausdorff series $\rho_1^\alpha(x)$:
\begin{equation*}
  \rho_1^\mu(x) = \mu - \sum_{k=1}^\infty \frac{1}{k!} \, \e_\mu^k(x) .
\end{equation*}
If $\g$ is a dg Lie algebra, only trees with vertices of valence $0$
or $1$ contribute to $\e_\alpha^k(x)$, and we recover the formula
\eqref{Deligne} figuring in the definition of the Deligne groupoid for
dg Lie algebras.

There is also a relative version of Theorem \ref{GAMMA}, analogous to
Theorem \ref{fibration}.
\begin{theorem} \label{Main}
  If $f:\g\to\h$ is a surjective morphism of nilpotent \Linf-algebras,
  the induced morphism
  \begin{equation*}
    \GAMMA_\bull(f):\GAMMA_\bull(\g)\to\GAMMA_\bull(\h)
  \end{equation*}
  is a fibration of simplicial sets.
\end{theorem}
\Proof
  Let $0\le i\le n$. Given a horn
  \begin{equation*} 
    \beta \in \ss(\Lambda^n_i,\GAMMA_\bull(\g))
  \end{equation*}
  and an $n$-simplex $\gamma\in\GAMMA_n(\h)$ such that
  \begin{equation*}
    f(\p_j\beta) = \p_j\gamma
  \end{equation*}
  for $j\ne i$, our task is to construct an element $\alpha\in
  f^{-1}(\gamma)\subset\GAMMA_n(\g)$ such that
  \begin{equation*}
    \p_j\alpha=\p_j\beta
  \end{equation*}
  if $j\ne i$. 
  
  Choose a solution $x\in\g^{1-n}$ of the equation $f(x)=I_{0\dots
    n}(\gamma)\in\h^{1-n}$. Let $\alpha$ be the unique element of
  $\GAMMA_n(\g)$ such that $\eps^i_n\alpha=\eps^i_n\beta$ and
  \begin{equation*}
    P_nR^i_n\alpha = (d+\delta) \Biggl( \sum_{k=1}^{n-1}
    \sum_{\substack{ i_1<\dots<i_k \\ i\notin\{i_1,\dotsc,i_k\} }}
    \om_{i_1\dotso i_k} \o I_{ii_1\dotso i_k}(\beta) + (-1)^i \,
    \om_{0\dotso\widehat{\imath}\dotso n} \o x \Biggr) .
  \end{equation*}
  If $j\ne i$, we have $\eps^i_n\p_j\alpha=\eps^i_n\p_j\beta$ and
  $P_nR^i_n\p_j\alpha=P_nR^i_n\p_j\beta$ and hence, by Lemma
  \ref{unique}, $\p_j\alpha=\p_j\beta$; thus, $\alpha$ fills the horn
  $\beta$. We also have
  $f(\eps^i_n\alpha)=f(\eps^i_n\beta)=\eps^i_n\gamma$ and
  $f(P_nR^i_n\alpha)=P_nR^i_n\gamma$, hence $f(\alpha)=\gamma$.
\Endproof

\begin{corollary}
  If $\g$ is a nilpotent \Linf-algebra, the inclusion of simplicial
  sets
  \begin{equation*}
    \GAMMA_\bull(\g) \hookrightarrow \MC_\bull(\g)
  \end{equation*}
  is a homotopy equivalence; in other words,
  $\pi_0(\GAMMA_\bull(\g))\cong\pi_0(\g)$, and for all $0$-simplices
  $\alpha_0\in\MC_0(\g)=\MC(\g)$,
  \begin{equation*}
    \pi_i(\GAMMA_\bull(\g),\alpha_0)\cong\pi_i(\g,\alpha_0) , \quad
    i>0 .
  \end{equation*}
\end{corollary}
\Proof
  This is proved by induction on the nilpotence length $\ell$ of $\g$.
  When $\g$ is abelian, $\MC_\bull(\g)$ and $\GAMMA_\bull(\g)$ are
  simplicial abelian groups, and their quotient is the simplicial
  abelian group
  \begin{equation*}
    \MC_n(\g)/\GAMMA_n(\g) \cong (d+\delta)s_n(\g\o\Om_n)^1 .
  \end{equation*}
  This simplicial abelian group is a retract of the contractible
  simplicial abelian group $\g\o\Om_\bull$, hence is itself
  contractible.
  
  Let $F^i\g$ be the lower central series of $\g$. Given $i>0$, we
  have a morphism of principal fibrations of simplicial sets
  \begin{equation*}
    \begin{CD}
      \GAMMA_\bull(F^{i+1}\g) @>>> \GAMMA_\bull(F^i\g) @>>>
      \GAMMA_\bull(F^i\g/F^{i+1}\g) \\
      @VVV @VVV @VVV \\
      \MC_\bull(F^{i+1}\g) @>>> \MC_\bull(F^i\g) @>>>
      \MC_\bull(F^i\g/F^{i+1}\g)
    \end{CD}
  \end{equation*}
  Since $F^i\g/F^{i+1}\g$ is abelian, we see that
  $\GAMMA_\bull(F^i\g/F^{i+1}\g)\simeq\MC_\bull(F^i\g/F^{i+1}\g)$. The
  result follows by induction on $\ell$.
\Endproof

When $\g$ is a nilpotent Lie algebra, the isomorphism
\begin{equation*}
  \pi_0(\GAMMA_\bull(\g)) \cong \pi_0(\MC_\bull(\g))
\end{equation*}
is equivalent to the surjectivity of the exponential map. The above
corollary may be viewed as a generalization of this fact.


\begin{thebibliography}{999}
  
\bibitem{Ashley} N. Ashley. \emph{Simplicial $T$-complexes and crossed
    complexes: a nonabelian version of a theorem of Dold and Kan.}
  Dissertationes Math. (Rozprawy Mat.) \textbf{265} (1988).
  
\bibitem{Beke} T. Beke. \emph{Higher \v{C}ech theory.} $K$-Theory
  \textbf{32} (2004), 293--322.
  
\bibitem{BG} A. K. Bousfield and V. K. A. M. Gugenheim. \emph{On
    $\textup{PL}$ de Rham theory and rational homotopy type.} Mem.
  Amer. Math. Soc. 8 (1976), no. 179.
  
\bibitem{Dakin} M. K. Dakin. ``Kan complexes and multiple groupoid
  structures.'' Esquisses Math., vol. 32, Univ. Amiens, Amiens, 1983.
  
\bibitem{Deligne} P. Deligne. Letter to L. Breen, February 28, 1994.
  Available (by permission of author) at
  \texttt{http://math.northwestern.edu/\textasciitilde
    getzler/Papers/deligne.pdf}.
  
\bibitem{Dold} A. Dold. \emph{Homology of symmetric products and other
    functors of complexes.} Ann. of Math. \textbf{68} (1958), 54--80.

\bibitem{Dupont} J. Dupont. \emph{Simplicial de Rham cohomology and
    characteristic classes of flat bundles.} Topology \textbf{15}
  (1976), 233--245.
  
\bibitem{DupontBook} J. Dupont. ``Curvature and characteristic
  classes.'' Lecture Notes in Math., no. 640, Springer-Verlag,
  Berlin-New York, 1978.
  
\bibitem{Duskin} J. Duskin. \emph{Higher-dimensional torsors and the
    cohomology of topoi: the abelian theory.} In ``Applications of
  sheaves.'' Lecture Notes in Math., no. 753, Springer, Berlin-New
  York, 1979, pp. 255--279.
  
\bibitem{DuskinNerve} J. Duskin. \emph{Simplicial matrices and the
    nerves of weak $n$-categories. I. Nerves of bicategories.} Theory
  Appl. Categ. \textbf{9} (2001), 198--308.
  
\bibitem{EM} S. Eilenberg and S. MacLane. \emph{On the groups
    $H(\Pi,n)$, I.} Ann. of Math. \textbf{58} (1953), 55--106.

\bibitem{FL} Z. Fiedorowicz and J.-L. Loday. \emph{Crossed simplicial
    groups and their associated homology.}  Trans. Amer. Math. Soc.
  \textbf{326} (1991), 57--87.

\bibitem{darboux} E. Getzler, \emph{A Darboux theorem for Hamiltonian
    operators in the formal calculus of variations.} Duke Math. J.
  \textbf{111} (2002), 535--560.

\bibitem{Glenn} P. G. Glenn. \emph{Realization of cohomology classes
    in arbitrary exact categories.} J. Pure Appl. Algebra 25 (1982),
  33--105.
  
\bibitem{GM} W. M. Goldman and J. J. Millson, \emph{The deformation theory
of representations of fundamental groups of compact K\"ahler manifolds.}
Inst. Hautes \'Etudes Sci. Publ. Math. \textbf{67} (1988), 43--96.

\bibitem{Hinich} V. Hinich. \emph{Descent of Deligne groupoids.}
  Internat. Math. Res. Notices (1997), 223--239.

\bibitem{Kan} D. M. Kan. \emph{Functors involving c.s.s. complexes.}
  Trans. Amer. Math. Soc. \textbf{87} (1958), 330--346.
  
\bibitem{LM} T. Lada and M. Markl. \emph{Strongly homotopy Lie
    algebras.} Comm. Algebra \textbf{23} (1995), 2147--2161.

\bibitem{LS} L. Lambe and J. Stasheff. \emph{Applications of
    perturbation theory to iterated fibrations.} Manuscripta Math.
  \textbf{58} (1987), 363--376.

\bibitem{May} J. P. May. ``Simplicial objects in algebraic topology.
  Reprint of the 1967 original.'' Chicago Lectures in Mathematics.
  University of Chicago Press, Chicago, IL, 1992.
  
\bibitem{NR} A. Nijenhuis and R. W. Richardson, Jr. \emph{Cohomology
    and deformations in graded Lie algebras.} Bull. Amer. Math. Soc.
  \textbf{72} (1966), 1--29.

\bibitem{Quillen} D. Quillen. \emph{Rational homotopy theory.} Ann. of
  Math. \textbf{90} (1969), 205--295.

\bibitem{Sullivan} D. Sullivan. \emph{Infinitesimal computations in
    topology.} Inst. Hautes \'Etudes Sci. Publ. Math. \textbf{47}
  (1977), 269--331.
  
\bibitem{Whitney} H. Whitney. ``Geometric integration theory.''
  Princeton University Press, Princeton, N.~J., 1957.

\end{thebibliography}
\end{document}